\documentclass[11pt]{article}
\usepackage[top=23.4mm, bottom=23.4mm, left=23.4mm, right=23.4mm]{geometry}
\usepackage{graphicx,dcolumn,bm}
\usepackage{amssymb,amstext,amsmath}
\usepackage{tikz-cd}
\usepackage{pifont}
\usepackage{graphicx}
\usepackage{dcolumn}
\usepackage{bm}
\usepackage{amssymb}
\usepackage{booktabs}  
\usepackage{siunitx}   
\usepackage{caption}   
\usepackage{array}     
\usepackage{multirow}
\usepackage{bigstrut}
\usepackage{makecell,rotating}
\usepackage{mathrsfs}
\usepackage{enumitem}
\usepackage{booktabs}
\usepackage{threeparttable}
\usepackage{multirow}
\usepackage{subfigure}
\usepackage{colortbl}
\usepackage{epsfig}
\usepackage{threeparttable}
\usepackage{chngpage}
\usepackage{float}
\usepackage{amstext}
\usepackage{amsmath}
\usepackage{amsthm}
\usepackage{pdflscape}
\usepackage{authblk}
\usepackage{hyperref}
\usepackage{appendix}
 \usepackage[utf8]{inputenc} 
\usepackage{rotating}
\usepackage[graphicx]{realboxes}
\usepackage{adjustbox}
\usepackage{setspace}
\usepackage{algorithm}
\usepackage{algorithmic}
\usepackage{caption}
\usepackage{xcolor}
\usepackage{xr}
\usepackage{caption}
\usepackage{lineno}
\allowdisplaybreaks
\newtheorem{lemma}{Lemma}
\newtheorem{corollary}{Corollary}
\newtheorem{theorem}{Theorem}
\newtheorem{assumption}{Assumption}
\newtheorem {Proposition}{Proposition}
\newtheorem{remark}{Remark}
\newtheorem{definition}{Definition}
\usepackage{enumitem}

\DeclareMathOperator{\rank}{rank}

\DeclareMathOperator{\trace}{Trace}
\DeclareMathOperator{\argmin}{argmin}
\DeclareMathOperator{\blkdiag}{blkdiag}
\DeclareMathOperator{\diag}{diag}
\DeclareMathOperator{\Sparse}{Sparse}
\DeclareMathOperator{\dist}{dist}
\DeclareMathOperator{\Struct}{Struct}
\DeclareMathOperator{\vvec}{vec}

\setitemize[1]{itemsep=1pt,partopsep=0pt,parsep=\parskip,topsep=4pt}
\title{Action-Dependent Policies for Multi-Agent Linear Quadratic Regulation}
 \author{Yongjie Liu$^{1}$, Gangshan Jing$^{2*}$, Long Wang$^{1*}$}
\date{
    $^1$ Center for Systems and Control, College of Engineering, Peking University, Beijing 100871, China\\[2ex]%
    $^2$ School of Automation, Chongqing University, Chongqing 400044, China\\[2ex]
    $^*$\footnotesize Corresponding authors. E-mail: jinggangshan@cqu.edu.cn,  longwang@pku.edu.cn\\[2ex]
}

\begin{document} 

\maketitle

\begin{abstract}
The linear quadratic regulation (LQR) problem over a multi-agent system usually requires a distributed feedback gain that is structurally constrained by the network topology, incurring an inevitable suboptimality gap. Existing studies either identify conditions on the system and topology that enable convex relaxation or seek a stationary point via policy gradient algorithms. In this paper, we reduce the suboptimality gap by adopting action-dependent policies (ADPs), where the control input of each agent depends on not only states of neighbor agents, but also control inputs of some other agents. We show that although exhibiting a distributed structure, ADP can enlarge the feasible set of feedback gains 
and even represent centralized optimal policies under suitable conditions. However, this comes at the price of losing coercivity of the cost function and possible failure of convergence for the projected gradient flow (PGF). To overcome this difficulty, we propose a partially regularized formulation, which guarantees convergence of PGFs, and forges a connection between standard multi-agent LQR and its action-dependent counterpart. Sufficient conditions for local convergence to globally optimal policies under PGFs regardless of regularization are also established.  Finally, numerical examples are presented to validate the performance of ADP.
\end{abstract}

\section{Introduction}
\label{sec:introduction}
In recent years, advances in policy optimization (PO) methods within reinforcement learning (RL) have renewed interest in linear quadratic regulation (LQR) problems \cite{hu2023}. For centralized LQR (no structural constraint on the control gain), the cost function is shown to be coercive and gradient dominant, ensuring global convergence of first-order methods \cite{fazel2018global,bu2019lqr,bu2020,poli2025}. Subsequent studies have extended these results to LQR with output feedback, risk-constraints, and model-free settings \cite{fatkhullin2021optimizing,zhao2023global,mohammadi2020linear,mohammadi2021convergence,jing2020model}. However, these results are limited to single-agent systems.

The widespread deployment of multi-agent systems (MASs), which typically consist of multiple coupled agents, has spurred growing attention to cooperative control and optimization problems \cite{jad2003, 2006flock,2007Saber,2009opt,Xiao2008,wang2010,2008xi}. As a canonical problem in this area, the multi-agent LQR problem aims to minimize a global quadratic cost for a network of coupled agents. When each agent has access only to local information, such as its own measurements and data received from its neighbors, the problem is referred to as \emph{distributed multi-agent LQR}, and the resulting controller is called a \emph{distributed controller}. Unlike its centralized counterpart, designing distributed optimal controllers is generally NP-hard \cite{cou1968}.  In particular, the authors in \cite{feng2020} showed that the set of stabilizing static decentralized controllers can possess exponentially many connected components. Despite these challenges arising from the lack of full state information, a number of notable works, including \cite{2002dis,rotkowitz2005characterization,nayyar2013decentralized,2015aut,8796100}, have appeared since the early 2000s. Existing approaches mainly fall into two categories: 1) identifying system characteristics and sparsity (or delay) conditions that admit tractable algorithms; 2) deriving suboptimal controllers with bounded suboptimality gaps. For the former, the quadratic invariance condition allows constrained minimum-norm problems to be solved via convex optimization \cite{rotkowitz2005characterization}.  For the latter, convex relaxation methods \cite{7464306,11264282} and sparsity invariance techniques \cite{8796100} have been developed. In short, much of the prior work involves a trade-off between the degree of decentralization and suboptimality \cite{shin2023,Zhang2025}. Nevertheless, a central question remains:

\textit{Given a multi-agent LQR problem, is it possible to achieve global optimality under a distributed control architecture?} 

Behind this question lies a \emph{policy representational bottleneck}: the centralized optimal gain is generally dense, whereas a distributed state feedback gain must respect the sparsity pattern of the sensing/communication graph. From this viewpoint, the difficulty essentially comes from the limited representation capability of the distributed policy class itself. The same bottleneck has been observed in multi-agent reinforcement learning (MARL). Intuitively, executing actions in a correlated manner leads to better cooperation; however, the commonly adopted independent/factorized policies cannot represent joint actions that require coordination.
 Until very recently, \emph{action-dependent policies} (ADP) were introduced in \cite{pmlr-v162-fu22d,pmlr-v202-chen23an,ding2025} to enlarge the joint policy space. Specifically, action dependencies refer to the dependence of an agent's decision-making on the actions of other agents. As shown in \cite{pmlr-v162-fu22d,pmlr-v202-chen23an,ding2025}, action dependencies provide a principled way to enlarge the representational capacity of distributed policies, and naturally lead to a better performance.

 Inspired by the aforementioned developments, we apply ADPs to the multi-agent LQR problem. In the literature, the usual distributed control input at each agent depends on its own state and the states of its neighbors through local sensing/communication. After introducing action dependencies,  the control input of one agent additionally depends on inputs of its neighbors in the \emph{action dependency graph} (ADG).
 Compared with densifying the sensing graph, the ADG offers extra design freedom, which naturally enlarges the policy space.
However, this enlargement does not come for free: the induced cost function is no longer coercive and the optimal control policies are generally non-isolated, so that  policy gradient (PG) methods may fail to converge. How to exploit the enlarged policy space while restoring convergence, and whether the enlarged space suffices to recover the centralized optimal policy, are topics addressed in this paper.

\subsection{Main contributions} 
\begin{itemize}[leftmargin=0pt]
\item \textbf{Action dependency graphs for fully dense policy representation:}  To characterize the enlarged policy space, the ADG is introduced to the multi-agent LQR problem for the first time, to the best of our knowledge.  We demonstrate the significant contribution of action dependencies to representing a fully dense policy. More precisely, the sparsity patterns of the induced gains are identified (see Proposition \ref{pro:spars}), and a necessary and sufficient graphical condition is established to characterize the representability for generic policies (see Theorem \ref{thm:c1}). 

\item \textbf{Convergence of policy gradient methods with action dependencies: }  Despite the enlargement of the policy space, we emphasize that projected gradient flow (PGF) for multi-agent LQR with action dependencies may not converge due to the lack of coercivity (see Example 2). Therefore, we introduce a partial regularization term to overcome this difficulty (see Theorem \ref{thm:reg}), and discover that the partially regularized problem acts as a novel bridge between the original multi-agent LQR and its action-dependent counterpart (see Theorem \ref{conne}). 

\item \textbf{Optimality guarantee and improvement of performance:} To achieve global optimality, we utilize the \emph{complementary sparsity conditions} on ADGs and sensing graphs to establish the guarantee in a generic sense (see Theorem \ref{thm:c3} and Theorem \ref{thm:c6}). The analyticity, \emph{Polyak–Łojasiewicz} inequality, and \emph{Morse-Bott} property of cost functions at stationary points are fully exploited to establish the exponentially asymptotic stability of the equilibrium set. Additionally,  an algorithm (see Algorithm \ref{alg:buildpar}) is proposed to guarantee the performance improvement over the PGF of the distributed multi-agent LQR problem. 
\end{itemize}
\subsection{Related work}
\begin{itemize}[leftmargin=0pt]
\item \textbf{Action-dependent policies in MARL:} Existing works concerning ADPs in MARL mainly include \cite{pmlr-v162-fu22d,ding2025,gupta2026}. All these works address the policy representational bottleneck and show that compared to independent
policies, ADPs enable more effective cooperation and even achieve global optimality under certain conditions.  However, these results rely crucially on a finite action space.
 Although the ADG adopted in this paper is similar to that in \cite{ding2025}, the continuous nature of the action space, and the complexity of the optimization landscape in the multi-agent LQR setting make their results not directly extendable to our case.

\item \textbf{Policies beyond state feedback in the multi-agent LQR:} Beyond the conventional state feedback, various policies, including output feedback  \cite{2020L4DC,rotkowitz2005characterization}, disturbance feedback \cite{2024ye,2015aut}, and  policies leveraging neighbors’ states and inputs \cite{zhang2020}, have
also been considered in the multi-agent LQR setting. In particular, \cite{zhang2020} designed distributed optimal controllers for the multi-agent LQR with each agent leveraging its neighbors' state and action information; yet they assume homogeneous dynamics and a consensus-oriented cost. In contrast, we consider general linear dynamics and establish sufficient conditions for achieving centralized optimality.

\item \textbf{Policy gradient algorithms in LQR:} 
Convergence results of PG methods for LQR  have been extensively
studied \cite{fazel2018global,bu2019lqr,fatkhullin2021optimizing,Dai2025,Guo2026,10494371}. In our work, the introduction of ADGs provides a reparameterization of controllers. Besides the loss of coercivity, redundancy in parameterization indicates that the optimal policies can form a smooth manifold. Recall that a similar non-isolation of optimizers arises in Linear Quadratic
Gaussian (LQG) problems, and is caused by the invariance of the cost under similarity transformations\cite{tang2023,quoti}; here it instead  stems from the invariance of the cost along the \emph{fiber} of the parameterization. Since the fiber is not well-behaved (its dimension may vary with the \emph{base point}), quotienting out symmetries as in \cite{quoti} is inapplicable, which distinguishes our analysis from prior work.

\end{itemize}
\subsection{Organization and notation}
\label{no:h}
\noindent\textbf{Organization:} The rest of this paper is organized as follows: Section 2 reviews preliminaries on the LQR problem, introduces ADPs and formulates the problem; Section 3 discusses policy representation; Section 4 investigates the optimization landscapes, and convergence behaviors for PGFs; Section 5 focuses on the globally optimal policy synthesis; Section 6 presents numerical simulations; Section 7 concludes the paper.

\noindent\textbf{Notation:} In this paper, $\mathbb{R}$ represents the set of real numbers;  $\mathbb{N}^{\ast}$ is the set of nonnegative integers; $I_{n}\in\mathbb{R}^{n\times n}$ is the unit matrix; $\mathbb{R}^{n}$ denotes the vector space of dimension $n$;  $\mathbf{1}_{n}\triangleq (1,\ldots,1)^{\top}\in\mathbb{R}^{n}$; For a set $T$, $|T|$ denotes the number of elements in $T$; Given two vector spaces $V_1$ and $V_2$, we write $V_1\cong V_2$
to indicate that they are linearly isometric;  The norm of $s=(s_1,...,s_n)^\top\in\mathbb{R}^{n}$ is denoted by $||s||\triangleq[\sum\limits_{k=1}^{n}{s_{k}^{2}}]^{\frac{1}{2}}$; For $X,Y\in\mathbb{R}^{m\times n}$, the Frobenius product on $\mathbb{R}^{m\times n}$ is defined by $\langle X,Y \rangle\triangleq \trace{(X^\top Y)}$, $X\leq (<) Y$ holds if and only if $X_{ij}\leq (<) Y_{ij}$ holds for $i=1,2,\dots,m$, $j=1,2,\dots,n$;  $\|X\|\triangleq \langle X,X \rangle ^{\frac{1}{2}}$
   denotes the Frobenius norm of $X$; For a matrix $X\in\mathbb{R}^{m\times n}$, the kernel space and rank of $X$ are denoted by  $\ker(X)$ and $\rank(X)$, respectively; $X^{\top}$ denotes the transpose of $X$;  $\vvec(X)$ is a column vector that stacks the columns of matrix $X$, and the inverse operator $\vvec^{-1}$ transforms a vector into a matrix with the suitable dimension; Given a space (or set) $\mathcal{K}$ of matrices, the vectorized space (or set) of $\mathcal{K}$ is $\vvec(\mathcal{K})\triangleq\{\vvec(K):K\in\mathcal{K}\}$; The sparsity structure of matrix $X$ is a Boolean matrix, denoted as $\Struct(X)$ with the same size of $X$ and $\Struct(X)_{ij}=1$ if $X_{ij}\neq 0$, otherwise, $\Struct(X)_{ij}=0$; 
 $X$ is called a fully dense matrix if all elements of $\Struct(X)$ equal $1$; The sparsity space associated with matrix $X\in\mathbb{R}^{m\times n}$ is $\Sparse(X)\triangleq\{Y\in \mathbb{R}^{m\times n}:Y_{ij}=0\ \text{if}\ X_{ij}=0\}$; Matrix $X$ satisfies $X>0$ indicates each element of $X$ is positive; Given $X\in\mathbb{R}^{m\times n}$, $\delta>0$, and subset $\mathcal{C}\subseteq\mathbb{R}^{m\times n}$, then $\mathcal{B}(X,\delta)\triangleq\{Y\in\mathbb{R}^{m\times n}:\|Y-X\|<\delta\}$ is called the $\delta$-neighborhood of $X$ in $\mathbb{R}^{m\times n}$, and  $\dist[X,\mathcal{C}]=\inf_{Y\in\mathcal{C}}\|X-Y\|$ is called the distance from $X$ to set $\mathcal{C}$;  For $X\in\mathbb{R}^{N\times N}$, its determinant is denoted by $\det(X)$, and '$X\succ (\succeq) 0$' means that $X$ is positive (semi-)definite; $\otimes$ is the Kronecker product;  For any positive integer $m$, let $H_{l,m}\in\mathbb{R}^{m\times m}$ denote the
lower triangular matrix whose entries on and below the main diagonal are all $1$
(remaining entries are $0$).

Consider a directed graph $\mathcal{G}=(\mathcal{V},\mathcal{E})$, where $\mathcal{V}=\left\{ {1,2,\ldots,N} \right\}$ is the vertex set, $\mathcal{E}\subseteq \mathcal{V}\times\mathcal{V}$ is the edge set. Let $\deg_{in}(i)$ be the in-degree of vertex $i$, and $\mathcal{N}_{i}\triangleq\{j\in\mathcal{V}:(j,i)\in\mathcal{E}\}$ be the set of neighbors of vertex $i$ in graph $\mathcal{G}$. A directed path is a sequence of edges $(v_{i_1},v_{i_2}),\ (v_{i_2},v_{i_3}),\dots,(v_{i_{k-1}},v_{i_k})$, with $(v_{i_{j-1}},v_{i_j})\in\mathcal{E}$ for $j=2,3,\dots,k$. $\mathcal{G}$ can be described by an adjacency matrix $(a_{ij})\in\mathbb{R}^{N\times N}$, where $a_{ij} = 1$ if $(j,i)\in\mathcal{E}$\footnote{Unlike the notation in \cite{FB}, we adopt the convention that  $a_{ij}=1$ indicates an edge from node $j$ to node $i$.}, but 0 otherwise. The transitive closure of $\mathcal{G}$  is the graph $\mathcal{G}^*=(\mathcal{V},\mathcal{E}^{\ast})$, where $(u,v)\in\mathcal{E}^{\ast}$ if there exists a directed path from $u$ to $v$ in $\mathcal{G}$.

\section{Preliminaries and problem formulation}
In this section, we will briefly review centralized and distributed multi-agent LQR problems, respectively.  Then, ADPs and their integration into multi-agent LQR will be introduced. 
\subsection{Centralized multi-agent LQR}
Here we focus on the MAS consisting of $N$ agents with dynamics:
\begin{equation}
    x_i(t+1)=\sum_{j=1}^{N}A_{ij}x_{j}(t)+B_i u_i(t),
\end{equation}
where $t\in\mathbb{N}^{\ast}$, $x_i\in\mathbb{R}^n$, $u_i\in\mathbb{R}^m$ ($m\leq n$), $A_{ij}\in\mathbb{R}^{n\times n}$, and $B_i\in\mathbb{R}^{n\times m}$ for $i,j=1,2,\dots,N$. The dynamics of the whole system can be expressed in a compact form: 
\begin{equation}\label{eqqq2}
    x(t+1)=Ax(t)+Bu(t),
\end{equation}
where $A\in\mathbb{R}^{nN\times nN}$ satisfies $A[(i-1)n+1:in,(j-1)n+1:jn]=A_{ij}$, $B=\blkdiag\{B_1,\dots,B_N\}\in\mathbb{R}^{nN\times mN}$, $x(t)=(x_1(t)^\top,\dots,x_N(t)^\top)^\top\in\mathbb{R}^{nN}$, and $u(t)=(u_1(t)^\top,\dots,u_N(t)^\top)^\top\in\mathbb{R}^{mN}$. 

Then we consider the following multi-agent LQR problem
\begin{align}
    \min\  \mathbb{E}\bigg[\sum_{t=0}^{\infty}x^{\top}(t)Qx(t)+u^{\top}(t)Ru(t)\bigg]\ 
    \text{s.t.}\  (\ref{eqqq2}),\ x(0)\sim \mathcal{D},
    \label{pro:lqrc}
\end{align}
where $Q\in\mathbb{R}^{nN\times nN}$ and $R=\blkdiag\{R_1,\dots,R_N\}\in\mathbb{R}^{mN\times mN}$ are positive definite, and $\mathcal{D}$ is a probability distribution of the initial state such that $x(0)$ is bounded and satisfies $\Sigma=\mathbb{E}[x(0)x^\top(0)]\succ0$.
In this paper, we first make the following assumption to make the LQR problem well-posed. 
\begin{assumption}\label{ass:sol}
    System $[A,B]$ is stabilizable.
\end{assumption}
Under Assumption \ref{ass:sol}, the centralized optimal solution to the LQR problem (\ref{pro:lqrc}) corresponds to a static linear state feedback controller 
$ u(t)=-Kx(t) $ ,
where $K=(R+B^{\top}PB)^{-1}B^{\top}PA$ is the optimal control policy, and $P$ is the unique positive semi-definite solution to the algebraic Riccati equation (ARE) \cite{lewis2012optimal} 
\begin{equation}\label{eq:are}
    P=A^\top PA+Q-A^\top PB (R+B^{\top}PB)^{-1}B^{\top}PA.
\end{equation}

From the perspective of PO, the LQR problem (\ref{pro:lqrc}) can be reformulated as \cite{fazel2018global}

\noindent\textbf{Problem $\mathbf{P}_0$:}
\begin{align}\label{pro:lq00}
    \min_{K} &\ J_0(K)=\trace(X\Sigma)\\
    \text{s.t.}\  & X=Q+K^\top RK+(A-BK)^\top X(A-BK),\  X\succ 0. \label{eq7}
\end{align}
Denote the set of stabilizing feedback gains $\mathcal{S}=\{K\in\mathbb{R}^{mN\times nN}: \rho(A-BK)<1\},$
where $\rho(\cdot)$ is the spectral radius of a matrix. In fact, $\mathcal{S}$ is the feasible set of $\mathbf{P}_0$, and is shown to be connected\cite{bu2019lqr}. Furthermore, the gradient dominance property always holds for $J_0$, which enables the convergence of PG methods.  For further results on PO, we refer interested readers to \cite{fazel2018global,bu2019lqr,bu2019}.

\subsection{Distributed multi-agent LQR}
 For an MAS with distributed control, the controller of each agent utilizes only local information from its neighbors. 
Given a graph $\mathcal{G}=(\mathcal{V},\mathcal{E})$ with $|\mathcal{V}|=N$ and $p,q\in\mathbb{N}^{*}$, define

\begin{equation}\label{defgs0}
\Sparse(\mathcal{G}^{p,q})\triangleq \Sparse(G\otimes \mathbf{1}_p\mathbf{1}^\top_q),
\end{equation}
 where $G$ is the adjacency matrix of graph $\mathcal{G}$.

For distributed controller synthesis, the \textit{sensing graph} is usually used to model the required sensing relationships \cite{10494371}.
\begin{definition}
   The sensing graph  $\mathcal{G}_S=(\mathcal{V},\mathcal{E}_S)$ is a directed graph with a self-loop at every node. An edge $(i,j)\in\mathcal{E}_S$ indicates that node $j$ has access to the state $x_i$ of node $i$. 
\end{definition}

Denote $\mathcal{K}=\Sparse(\mathcal{G}_{S}^{m,n}).$ Then the distributed LQR problem can be formulated as 

\noindent \textbf{Problem $\mathbf{P}_1$:}
\begin{equation}\label{pro:lq0}
    \begin{aligned}
    \min_{K} \ J_d(K)=\trace(X\Sigma),
    \text{s.t.}\  (\ref{eq7}),\ K\in\mathcal{K}.
\end{aligned}
\end{equation}

To ensure the feasibility of Problem $\mathbf{P}_1$, we propose the following assumption: 
\begin{assumption}\label{ass:st}
    $\mathcal{K}\cap\mathcal{S}\neq\emptyset.$
\end{assumption}

Notably, different from Problem $\mathbf{P}_0$, the feasible set $\mathcal{K}\cap\mathcal{S}$ can be disconnected\cite{feng2020}. This implies that without prior knowledge for initialization, local search methods generally cannot solve $\mathbf{P}_1$ with global optimality guarantee\cite{feng2020}. Further results on PO can be found in \cite{bu2019lqr,feng2020}.

\subsection{Action-dependent policies}
\label{prob}


In this paper, we consider multi-agent LQR problems with a control policy for each agent $i$ in the following form:

\begin{equation}\label{eq:con1}  
    u_{i}(t)=-\sum\limits_{(j,i)\in\mathcal{E}_S}K_{ij}x_j(t)+\sum\limits_{(j,i)\in\mathcal{E}_A}L_{ij}u_j(t),
\end{equation}
where the second part on the right hand side takes the action dependencies into account, $K_{ij}\in\mathbb{R}^{m\times n}$ is the feedback gain using  the state of agent $j$,  $L_{ij}\in\mathbb{R}^{m\times m}$ is the gain using the action of agent $j$, $\mathcal{E}_S$ represents the edge set of the sensing graph, and $(j,i)\in\mathcal{E}_A$ indicates that agent $i$ has access to the input $u_j$ of agent $j$. In practice, the transmission of the action information between two agents can be realized by communications.  

Note that action dependency relationships are described by the \emph{action dependency graph} (ADG) \cite{ding2025} in the following definition. In order to avoid deadlocks in decision-making of MASs, it is natural to assume that ADGs are directed and acyclic. %
\begin{definition}
   The action dependency graph $\mathcal{G}_A=(\mathcal{V},\mathcal{E}_A)$ is a directed acyclic graph in which $(i,j)\in\mathcal{E}_A$ indicates that node $j$ has access to the action $u_i$ of node $i$. 
\end{definition}

Due to the directed acyclic property of $\mathcal{G}_A$, there exists at least one topological sort such that agents perform their actions in the corresponding order \cite{ding2025}. Without loss of generality, we adopt a labeling scheme of the vertex set such that for each $(i,j)\in\mathcal{E}_A$, it holds that $i<j$\footnote{{A different relabeling of nodes corresponds to a permutation matrix $W\in\mathbb{R}^{N\times N}$ satisfying $WW^\top=I_N$ (see \cite[Remark 4.6]{FB}). The adjacency matrices of $\mathcal{G}_A$ and $\mathcal{G}_S$ become $WG_AW^\top$ and $WG_SW^\top$, respectively.}}. Therefore, the adjacency matrix of $\mathcal{G}_A$ is lower triangular with zero diagonal entries. 
Denote $\ \mathcal{L}=\Sparse(\mathcal{G}_{A}^{m,m})$. Based on (\ref{defgs0}),  the compact form of (\ref{eq:con1}) is 
\begin{equation}
u=-Kx+Lu,\ K\in\mathcal{K},\ L\in\mathcal{L}. 
\end{equation}Due to the specific labeling of the vertex set of $\mathcal{G}_A$, $L$ is a lower triangular matrix, and its diagonal blocks are zero, which indicates that matrix $I-L$ is always non-singular. Hence Equation (\ref{eq:con1}) can be equivalently written as
\begin{equation}
    u=-(I-L)^{-1}Kx,\ K\in\mathcal{K},\ L\in\mathcal{L}.
\end{equation} 
Denote the space of control gains applied to the system state directly by: 
\begin{equation}
    \mathcal{M}=\{(I-L)^{-1}K: K\in\mathcal{K}, L\in\mathcal{L}\}. 
\end{equation}
Consider mapping $T: \mathcal{L}\times\mathcal{K}\to\mathcal{M}\subseteq\mathbb{R}^{mN\times nN}$ defined by 
\begin{equation}\label{deft}
    T(L,K)=(I-L)^{-1}K.
\end{equation}

Now we are ready to consider the following LQR problem with ADPs.

\noindent\textbf{Problem $\mathbf{P}_2$:}
\begin{equation}\label{pro:lqadg}
    \begin{aligned}
    \min_{K,L} &\ J(L,K)=\trace(X\Sigma)\\
    \text{s.t.}\  & X=Q+\tilde{K}^\top R\tilde{K}+(A-B\tilde{K})^\top X(A-B\tilde{K}),\  X\succ 0,\\
    & \tilde{K}=T(L,K), K\in\mathcal{K}, L\in\mathcal{L}.
\end{aligned}
\end{equation}
 Denote the feasible set of $\mathbf{P}_2$ as $$\mathcal{F}\triangleq\{(L,K)\in\mathcal{L}\times\mathcal{K}:\rho(A-B(I-L)^{-1}K)<1\}.$$ Note that if $\mathcal{L}=0$,  this problem reduces to Problem $\mathbf{P}_1$. Thus, $\mathcal{K}=T(0,\mathcal{K})\subseteq\mathcal{M}$ always holds. Then Assumption \ref{ass:st} indicates $\mathcal{M}\cap\mathcal{S}\neq\emptyset$ and $\mathcal{F}\neq\emptyset$.  We point out that  $\mathcal{F}$ provides a larger search space for structural controller synthesis, potentially enabling a better performance.

 For each $\tilde{K}\in\mathcal{M}$, $T^{-1}(\tilde{K})$ is called a \textit{fiber} at $\tilde{K}$, and $\tilde{K}$ is called the \emph{base point} of this fiber. Note that $J$ is constant on each fiber, i.e., for any $(L,K)\in T^{-1}(\tilde{K})$, we have
\begin{equation}\label{eq:inv}
J(L,K)=J_0(\tilde{K}).
\end{equation}
More detailed properties of the fiber will be discussed in the next section. Specifically, we will show that the introduction of action dependencies enlarges the policy space at the expense of breaking the coercivity of the cost function. 
Therefore, to handle the lack of coercivity, we introduce the following partially regularized problem with a penalty parameter $\lambda>0$.

\noindent\textbf{Problem $\mathbf{P}_3$:}
\begin{equation}\label{pro:lqadg}
    \begin{aligned}
    \min_{K,L} &\ J_\lambda(L,K)=\trace(X\Sigma)+\frac{\lambda}{2}\|L\|^2\\
    \text{s.t.}\  & X=Q+\tilde{K}^\top R\tilde{K}+(A-B\tilde{K})^\top X(A-B\tilde{K}),\ X\succ 0, \\
    & \tilde{K}=T(L,K), K\in\mathcal{K}, L\in\mathcal{L}.
\end{aligned}
\end{equation}
\begin{figure}[htbp!]
\centering
\includegraphics[width=0.7\textwidth]{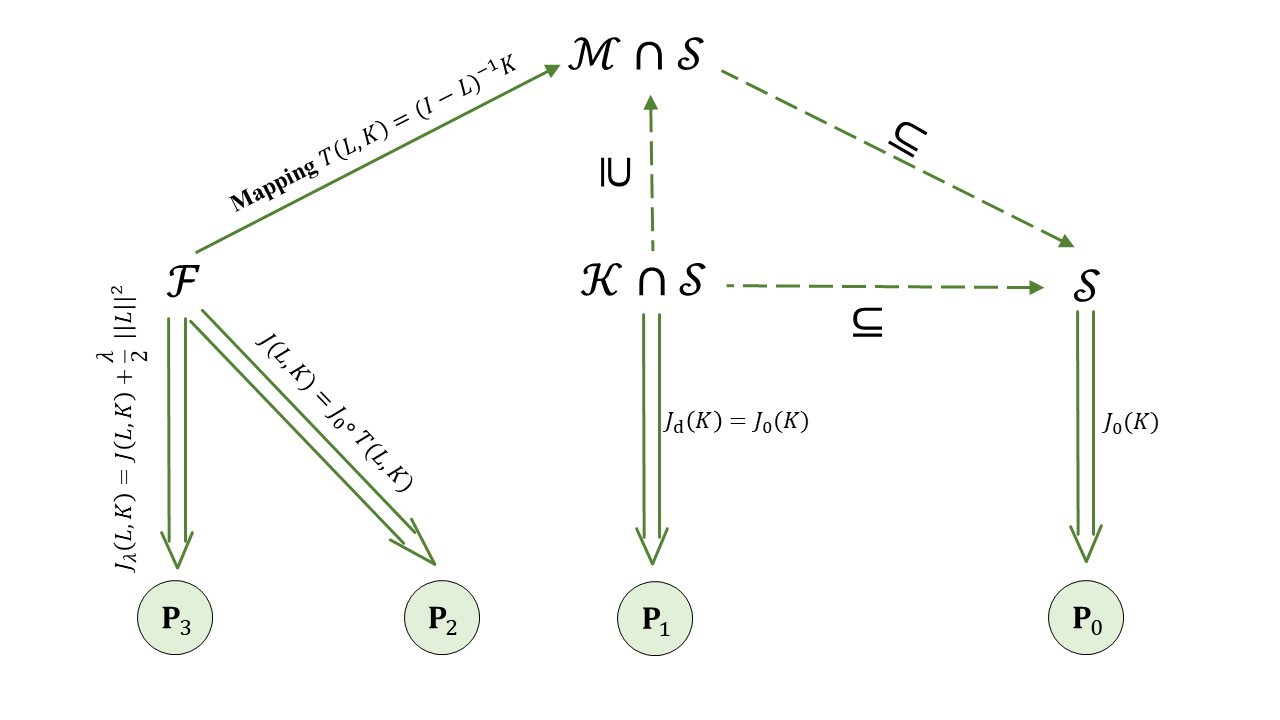}
\caption{Relationships among four problems.  The upper part shows the relations among these feasible domains through the inclusion  and the mapping $T$, while the lower part presents the relations among these cost functions. Note that Problems $\mathbf{P}_2$ and $\mathbf{P}_3$ share the same feasible domain $\mathcal{F}$, while the feasible domains of Problems $\mathbf{P}_0$ and $\mathbf{P}_1$ are $\mathcal{S}$ and $\mathcal{K}\cap\mathcal{S}$, respectively. } 
\label{fig:grr1}
\end{figure}

Up to now, we have introduced Problems $\mathbf{P}_0$-$\mathbf{P}_3$. Figure 1 presents a diagram to illustrate the relationships among the proposed problems and their feasible sets. In summary, we will concentrate on the following question throughout this work.

\noindent\textbf{Question:} \textit{Compared with Problem $\mathbf{P}_1$, do ADPs contribute to performance improvement of $\mathbf{P}_2$ and $\mathbf{P}_3$? If yes, is it possible for an ADP to represent a globally optimal controller of $\mathbf{P}_0$ through suitable design of ADGs?}


\section{Policy representation under different ADGs }
According to Subsection \ref{prob},  $\mathcal{M}$ exactly represents policies incorporating action dependencies. In this section, we first identify the sparsity pattern of $\mathcal{M}$ through graph structures of $\mathcal{G}_S$ and $\mathcal{G}_A$. Then we discuss the policy representation capability under different ADGs.

\subsection{Sparsity structure of $\mathcal{M}$}
As shown in Subsection \ref{prob}, adjacency matrix $G_A$ is lower triangular with zero diagonal entries. Thus, the following nilpotency of $L\in \mathcal{L}$ holds naturally.
\begin{lemma}\label{lem:spl}
    If $L\in\mathcal{L}$,  then $L^{N}=0$, and 
\begin{equation}\label{eq:inv}
       (I-L)^{-1}=\sum\limits_{k=0}^{N-1}L^k.
\end{equation}
\end{lemma}
\begin{proof}
    Since $\mathcal{L}=\Sparse(G_A\otimes\mathbf{1}_m\mathbf{1}^\top_m)$, $L$ is a block lower triangular matrix with zero diagonal blocks, which indicates that $L^{N}=0$. Hence (\ref{eq:inv}) follows from  $(I-L)\sum\limits_{k=0}^{N-1}L^k=I-L^N=I$.
\end{proof}

Then we recall a useful lemma characterizing the sparsity pattern of the matrix product, which can be easily verified from the definition of matrix multiplication\cite{8796100}.
\begin{lemma}\label{lem:sp}
    Suppose $A_1\in\mathbb{R}^{p\times q}$ and $A_2\in \mathbb{R}^{q\times r}$ are Boolean matrices, if $X_i\in\Sparse(A_i)$ for $i=1,2$, then $X_1 X_2\in\Sparse(A_1 A_2)$.
\end{lemma}

Equivalently, Lemma \ref{lem:sp} shows that $\Struct(X_i)\leq A_i$ for $i=1,2$ implies $\Struct(X_1 X_2)\leq A_1 A_2$. Based on Lemma \ref{lem:spl} and Lemma \ref{lem:sp}, we can deduce the sparsity pattern of $\mathcal{M}$.
\begin{Proposition}\label{pro:spars}
    The following statements hold for $\mathcal{M}$:

       1) $\mathcal{M}$ is connected;
       
       2) Given $\tilde{K}\in\mathcal{M}$, $t\tilde{K}\in\mathcal{M}$ for $t\in\mathbb{R}$, and fiber $T^{-1}(\tilde{K})$ is an affine subspace of $\mathcal{L}\times\mathcal{K}$;

       3) $\mathcal{M}\subseteq \Sparse\Bigg((\sum\limits_{k=0}^{N-1} G^k_A G_S)\otimes \mathbf{1_m 1^\top_n}\Bigg)$, and there exists a matrix $K\in \mathcal{M}$ such that $\Struct(K)=\Struct\Bigg((\sum\limits_{k=0}^{N-1}G^k_A G_S)\otimes \mathbf{1_m 1^\top_n}\Bigg)$.
   \label{prop:sp}
\end{Proposition}
\begin{proof}
    See Appendix \ref{pro1}.
\end{proof}

Consequently, we identify three special cases: (1) Graph $\mathcal{G}_A$
   makes no contribution to policies in set $\mathcal{M}$ beyond those in $\mathcal{K}$; (2) Set $\mathcal{M}$ contains a fully dense matrix; (3) Set $\mathcal{M}$ exactly equals $\mathbb{R}^{mN\times nN}$. These cases are summarized in the following corollary.

\begin{corollary}\label{cor1}
    1) $\mathcal{M}=\mathcal{K}$  if and only if $\Struct(G_A G_S )\leq G_S$;

    2) $\mathcal{M}$ contains a fully dense matrix if and only if $\sum\limits_{k=0}^{N-1}G^{k}_A G_S >0$;

    3) $\mathcal{M}= \mathbb{R}^{mN\times nN}$ if and only if $\mathcal{G}_S$ is a complete graph.
\end{corollary} 
\begin{proof}
    See Appendix \ref{pro1}.
\end{proof}

Based on Corollary \ref{cor1}, we know that condition $\sum\limits_{k=0}^{N-1}G^{k}_A G_S>0$ is both necessary and sufficient  for $\mathcal{M}$ to contain a fully dense policy. For convenience, we give the following definition.
\begin{definition}
    ADG $\mathcal{G}_A$ and the sensing graph $\mathcal{G}_S$ are said to satisfy complementary sparsity condition I if  $\sum\limits_{k=0}^{N-1}G^{k}_A G_S>0$.
\end{definition}

In fact, 
 $\Struct(\sum\limits_{k=0}^{N-1}G^{k}_A)$ is the adjacency matrix of the transitive closure of $\mathcal{G}_A$. Intuitively, complementary sparsity condition I means that if $\mathcal{G}_S$ is sparse,  $\mathcal{G}_A$ should be dense to ensure the existence of a fully dense policy. This explains the name ``complementary sparsity". 
Notably, this condition does not ensure that $\mathcal{M}$ contains all fully dense policies. To illustrate this, we provide the following example. 

\begin{figure}[htbp!]
\centering
\includegraphics[width=0.8\textwidth]{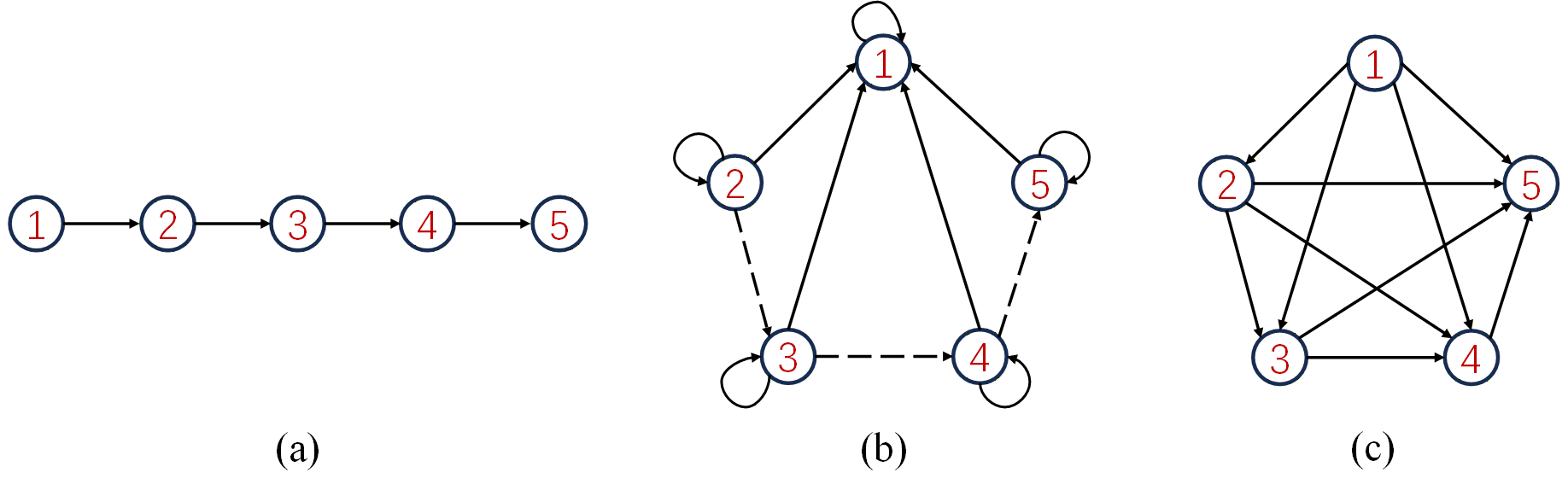}
\caption{Illustration of complementary sparsity condition I. (a) ADG $\mathcal{G}_A$ corresponding to a directed tree. Note that node $1$ is the root node of the tree, which indicates that input information of other nodes is not available to node $1$. (b) Two sensing graphs $\mathcal{G}_S$ and $\mathcal{G}'_S$, where edges of $\mathcal{G}_S$ correspond to the solid lines, while those of $\mathcal{G}'_S$ additionally include the dashed lines. In both cases, complementary sparsity condition I is valid here.  To guarantee $ H_{l,5} G_S>0$, all state information of other nodes should be available to node $1$ in the sensing graph. This reflects the sparse compensation characteristics of graphs $\mathcal{G}_A$ and $\mathcal{G}_S$. (c) The transitive closure of $\mathcal{G}_A$.}  
\label{fig:gr1}
\end{figure}

\noindent\textbf{Example 1:} Given adjacency matrices of ADG $\mathcal{G}_A$ and two sensing graphs corresponding to Figure \ref{fig:gr1} 
\begin{align}\label{eq:gas}
 G_A&= \begin{bmatrix}
0 & 0 & 0 & 0 & 0 \\
1 & 0 & 0 & 0 & 0 \\
0 & 1 & 0 & 0 & 0 \\
0 & 0 & 1 & 0 & 0 \\
0 & 0 & 0 & 1 & 0 
\end{bmatrix},\ G_S= \begin{bmatrix}
1 & 1 & 1 & 1 & 1 \\
0 & 1 & 0 & 0 & 0 \\
0 & 0 & 1 & 0 & 0 \\
0 & 0 & 0 & 1 & 0 \\
0 & 0 & 0 & 0 & 1
\end{bmatrix},\ 
G'_S = \begin{bmatrix}
1 & 1 & 1 & 1 & 1 \\
0 & 1 & 0 & 0 & 0 \\
0 & 1 & 1 & 0 & 0 \\
0 & 0 & 1 & 1 & 0 \\
0 & 0 & 0 & 1 & 1
\end{bmatrix}.
\end{align}
Take $L_0=G_A$, $K_0=G_S$, and $K'_0=G'_S$. Then a direct calculation yields that $X_0=(I-L_0)^{-1}K_0>0$, $X'_0=(I-L_0)^{-1}K'_0>0$. Specifically,
\begin{equation}\label{eq:gas1}
X_0= \begin{bmatrix}
     1  &   1  &   1  &   1   &  1\\
     1  &   2   &  1     &1   &  1\\
     1   &  2 &    2 &    1&     1\\
     1    & 2     &2    & 2 &    1\\
     1    & 2     &2     &2&     2
\end{bmatrix},\ X'_0= \begin{bmatrix}
     1  &   1  &   1  &   1   &  1\\
     1  &   2   &  1     &1   &  1\\
     1   &  3 &    2 &    1&     1\\
     1    & 3     &3    & 2 &    1\\
     1    & 3     &3     &3&     2
\end{bmatrix}.
\end{equation}
 Moreover, one can verify that $\sum_{k=0}^{4}G^k_A=H_{l,5}$, which corresponds to the transitive closure of $\mathcal{G}_A$ in Figure  \ref{fig:gr1}(c), and complementary sparsity condition I holds here. This is consistent with Corollary \ref{cor1}.  However, policy $X'_0$ can not be represented by $(I-L)^{-1}K$ with $(L,K)\in \Sparse(\mathcal{G}_{A})\times \Sparse(\mathcal{G}_{S})$. 
This can be checked by solving the following convex optimization problem:
\begin{equation}\label{pro:dec0}
    \min_{K,L} \ \|K-(I-L)X'_0\|^2 \ \
    \text{s.t.}\   L\in \Sparse(\mathcal{G}_{A}),K\in\Sparse(\mathcal{G}_{S}).
\end{equation}
 By using the CVXPY toolbox in Python, the minimal value for (\ref{pro:dec0}) is shown to be 0.9607, which indicates that $X'_0\notin\mathcal{M}$. 
\hfill$\square$


\begin{remark}
 Note that $\Struct(\sum\limits_{k=0}^{N-1}G^k_A)\leq H_{l,N}$, where $H_{l,N}$ is defined in Subsection \ref{no:h}.
Hence complementary sparsity condition I implies  $H_{l,N}  G_S >0$.
Consequently, the first row of $G_S$
   must be positive.  On the other hand, if the underlying undirected graph of the transitive closure of $\mathcal{G}_A$ is a complete graph, we have $\Struct(\sum\limits_{k=0}^{N-1}G^k_A)=H_{l,N}$, and then $H_{l,N}  G_S >0$ also implies complementary sparsity condition I. 
\end{remark}

\subsection{Policy representation}
    Inspired by Proposition \ref{prop:sp},  we seek a practical criterion for determining the representability of a policy $\tilde{K}$ through action dependencies, i.e., $\tilde{K}\in\mathcal{M}$. In fact, the key lies in solving the following matrix equation with structural constraints:
    \begin{equation}\label{eq:pa}
        (I-L)^{-1}K=\tilde{K}, (L,K)\in\mathcal{L}\times\mathcal{K}.
    \end{equation} 


To analyze \eqref{eq:pa}, consider the vectorized subspaces
$\vvec(\mathcal{K})\subseteq\mathbb{R}^{mnN^2}$
and $\vvec(\mathcal{L})\subseteq\mathbb{R}^{m^2N^2}$, and let
$\{\hat{K}_j\}_{j=1}^{n_k}$ and $\{\hat{L}_j\}_{j=1}^{n_l}$ be their bases, respectively.
Define
$\hat{K}=[\hat{K}_1,\dots,\hat{K}_{n_k}]\in\mathbb{R}^{mnN^2\times n_k}$
and $\hat{L}=[\hat{L}_1,\dots,\hat{L}_{n_l}]\in\mathbb{R}^{m^2N^2\times n_l}$.


\begin{Proposition}[Representability of a specific policy]\label{pron}
    Given $\tilde{K}\in\mathbb{R}^{mN\times nN}$, $\tilde{K}\in\mathcal{M}$  if and only if  
    \begin{equation}\label{rankc}
        \rank \bigg[\hat{K},(\tilde{K}^\top\otimes I_{mN})\hat{L}\bigg]=\rank \bigg[\hat{K},(\tilde{K}^\top\otimes I_{mN})\hat{L},\vvec(\tilde{K})\bigg].
    \end{equation}
Furthermore, given any $\tilde{K}\in\mathcal{M}$, the following holds
\begin{equation}
\dim(T^{-1}(\tilde{K}))=n_l+n_k-\rank \bigg[\hat{K},(\tilde{K}^\top\otimes I_{mN})\hat{L}\bigg].    
\end{equation}
\end{Proposition}
\begin{proof}
See Appendix \ref{appdis0}.
\end{proof}

   Notably, according to Corollary \ref{cor1},  $\mathcal{M}\subsetneqq\mathbb{R}^{mN\times nN}$ if $\mathcal{G}_S$ is not a complete graph. Naturally, how large can $\mathcal{M}$ be? To answer this question, we propose the following condition.

\begin{assumption}[Complementary sparsity condition II]\label{ass:den}
    $m\deg_{in,A}(i)+n\deg_{in, S}(i)\geq nN, \forall\ i=1,\dots,N.$
\end{assumption}


Now we show that Assumption \ref{ass:den} is necessary and sufficient for  $\mathcal{M}$ to contain almost every policy in $\mathbb{R}^{mN\times nN}$. 
For simplicity, let $A_i=\mathcal{N}_{in,A}(i)$, $S_i=\mathcal{N}_{in,S}(i)$, and  $S^c_i=\mathcal{V}\setminus S_i$. Denote $\delta_i=N-|S_i|$, $d_i=|A_i|$.  Given $\tilde{K}\in\mathbb{R}^{mN\times nN}$, we let $\tilde{K}_i=\tilde{K}\big((i-1)m+1:im,:\big)\in \mathbb{R}^{m\times nN}$, and define linear maps $\mathcal{Y}_i:\mathbb{R}^{mN\times nN}\to \mathbb{R}^{md_i\times n\delta_i}$ and $\mathcal{Z}_i:\mathbb{R}^{mN\times nN}\to \mathbb{R}^{m(d_i+1)\times n\delta_i}$ by
\begin{equation} \label{defyz}
    \mathcal{Y}_i(\tilde{K})\triangleq \bigg[(\tilde{K}_{j_1}\Big|_{S^c_i})^\top,\dots,(\tilde{K}_{j_{d_{i}}}\Big|_{S^c_i})^\top\bigg]^\top,\ 
    \mathcal{Z}_i(\tilde{K})\triangleq \begin{bmatrix}\mathcal{Y}_i(\tilde{K})\\ \tilde{K}_{i}\Big|_{S^c_{i}}\end{bmatrix},
\end{equation}
respectively, where $\{j_1,\dots,j_{d_{i}}\}=A_i$, and $\tilde{K}_{\cdot}\Big|_{S^c_{i}}$ represents the submatrix of $\tilde{K}_{\cdot}$ whose column blocks correspond to nodes in $S^c_i$. Furthermore, we introduce 
$S_0\subseteq\mathbb{R}^{mN\times nN}$ defined by
 \begin{align}\label{defs0}
       \ S_0\triangleq\bigcup_{i=1}^{N}S^i_{0},\
       S^i_{0}= \{\tilde{K}\in\mathbb{R}^{mN\times nN}:  \rank \big(\mathcal{Y}_i(\tilde{K})\big)<n \delta_i\}.
    \end{align}
Under Assumption~\ref{ass:den}, we will show that $S_0$  has (Lebesgue) measure zero.
Consequently, if $\tilde{K}\in\mathbb{R}^{mN\times nN}$ is generic\footnote{For a configuration or vector, genericity means algebraic independence of its coordinates\cite[Section 7.2]{lyjbb9.14}. After vectorizing matrix $\tilde{K}$, we obtain that genericity of $\tilde{K}$ implies $\tilde{K}\in \mathbb{R}^{mN\times nN}\setminus S_0$. In other words, the condition $\tilde{K}\in \mathbb{R}^{mN\times nN}\setminus S_0$ is weaker than the usual definition of genericity.}, then
$\tilde{K}\in\mathbb{R}^{mN\times nN}\setminus S_0$.
\begin{theorem}[Representability of generic policies]  \label{thm:c1}The following statements hold for $\mathcal{M}$: 

1) $\mathcal{M} = \bigl\{ \tilde{K} \in \mathbb{R}^{mN \times nN} :
  \rank\bigl( \mathcal{Y}_i(\tilde{K}) \bigr)
  = \rank\bigl( \mathcal{Z}_i(\tilde{K}) \bigr),\  \forall\, i = 1, \dots, N \bigr\}$.
  
2) $\mathbb{R}^{mN \times nN} \setminus S_0 \subseteq \mathcal{M}$.
Furthermore, $\mathbb{R}^{mN\times nN}\setminus S_0$ is dense and open in $\mathbb{R}^{mN\times nN}$ if and only if Assumption \ref{ass:den} holds.  
\end{theorem}

\begin{proof}
See Appendix \ref{appdis0}.
\end{proof}
 
\begin{corollary} \label{corr}
     $\mathcal{M}$ contains all policies in $\mathbb{R}^{mN\times nN}$
 except for a set of measure zero if and only if Assumption~\ref{ass:den} holds.
\end{corollary}
\begin{proof}
See Appendix \ref{appdis0}.
\end{proof}

Next, we discuss whether set $\mathcal{M}$ admits the structure of a \emph{smooth manifold}, which is helpful for later discussion on PO.  Generally, the answer is negative. 
Throughout the following analysis, we assume that Assumption \ref{ass:den} is valid.
Consequently, for any $\tilde{K}\in \mathbb{R}^{mN\times nN}\setminus S_0$, there exists $\delta>0$ such that $\mathcal{B}(\tilde{K},\delta)\subseteq\mathcal{M}$.  Then Proposition \ref{prop:sp} indicates that $t \mathcal{B}(\tilde{K},\delta)\subseteq \mathcal{M}$ holds for $t\in\mathbb{R}$. Hence, $\mathcal{M}$ contains the double cone $\mathcal{C}\triangleq\bigcup_{t\in\mathbb{R}}t\,\mathcal{B}(\tilde{K},\delta)$,
whose apex is the origin; in particular, $0\in\mathcal{M}$ and
$\mathcal{M}$ has nonempty interior in $\mathbb{R}^{mN\times nN}$.  However, based on Corollary \ref{cor1}, there exists $\tilde{K}^{\prime}\in\mathbb{R}^{mN\times nN}$ such that $t\tilde{K}^{\prime}\notin\mathcal{M}$ for $t\in\mathbb{R}\setminus\{0\}$ if $\mathcal{G}_S$ is not a complete graph. Therefore, $0$ is always a bad point such that there does not exist a neighborhood $\mathcal{U}_0$ of $0$ which is homeomorphic to an open subset of $\mathbb{R}^{mN\times nN}$.  On the other hand, Theorem \ref{thm:c1} also brings good news, i.e., $\mathcal{M}$ admits a local manifold structure near $\tilde{K}\in\mathbb{R}^{mN\times nN}\setminus S_0$ since $\mathcal{B}(\tilde{K},\delta)\cap\mathcal{M}=\mathcal{B}(\tilde{K},\delta)$ holds for some $\delta>0$. Therefore, we obtain the following fact. 

\noindent\textbf{Fact:}\label{cor2} \emph{In general, $\mathcal{M}$ fails to be a manifold globally;  under Assumption \ref{ass:den}, $\mathcal{M}$ admits a local manifold structure near $\tilde{K}\in\mathbb{R}^{mN\times nN}\setminus S_0$.}



\begin{remark} \label{re:ass} Concerning Theorem \ref{thm:c1}, the following comments are provided.

1) Since  $m \leq n$ typically holds in practice, Assumption \ref{ass:den} amounts to requiring that the sum of in-degrees for each node $i$ in graphs $\mathcal{G}_A$ and $\mathcal{G}_S$ is always larger than $N$. In particular, the converse also holds if $m=n$.  
    
2) Under Assumption \ref{ass:den}, Theorem \ref{thm:c1} implies that there exist fully dense matrices in $\mathcal{M}$. When $m=n=1$,  it follows from Corollary \ref{cor1} that $\sum\limits_{k=0}^{N-1}G^{k}_A G_S>0$. Therefore, complementary sparsity condition II implies complementary sparsity condition I, but not vice versa (see Figure \ref{fig:gr1}). 

3) Note that the dimension of each fiber $\dim(T^{-1}(\tilde{K}))$ can vary as the base point $\tilde{K}\in\mathcal{M}$ varies, which indicates that fibers are not linearly isometric in general. For instance, based on Proposition \ref{pron}, $\dim(T^{-1}(0))=n_l\leq\frac{N(N-1)}{2}m^2$. However, under Assumption \ref{ass:den}, we have $\dim(T^{-1}(\tilde{K}))\geq mn N^2-\rank \bigg[\hat{K},(\tilde{K}^\top\otimes I_{mN})\hat{L}\bigg]$.
\end{remark}



\section{Model-based policy optimization}
In this section, we employ PO for seeking the optimal action-dependent policy of the multi-agent LQR problem. First, the optimization landscapes of $\mathbf{P}_2$ and $\mathbf{P}_3$ will be investigated. Then, convergence behaviors of the PGFs will be discussed. The advantage of the partially regularized model $\mathbf{P}_3$ will be illustrated through an interesting example.
\subsection{Optimization landscape under action-dependencies}
Compared with Problem $\mathbf{P}_1$, Problems $\mathbf{P}_2$ and $\mathbf{P}_3$ share the same feasible domain $\mathcal{F}$. Furthermore, we will show that the feasible domain $\mathcal{F}$ can still be disconnected, and coercivity of cost function $J$ breaks down. 
\begin{Proposition}[Disconnectivity]\label{pro:dis}
 The feasible set $\mathcal{F}$ is disconnected.
\end{Proposition}
\begin{proof}
See Appendix \ref{appdis1}.
\end{proof}

    \begin{Proposition}[Coercivity]
 The cost function $J$ is not coercive on $\mathcal{F}$, while $J_\lambda$ is  coercive on $\mathcal{F}$.
    \label{pro:no}
\end{Proposition}
\begin{proof}
See Appendix \ref{appdis1.1}.
\end{proof}

As a significant difference with $\mathbf{P}_1$, the lack of coercivity in $\mathbf{P}_2$ suggests that $J(L,K)$ does not necessarily approach infinity as $(L,K)$ becomes unbounded. Based on Proposition \ref{prop:sp}, for $\tilde{K}\in\mathcal{M}$, $T^{-1}(\tilde{K})$ is either a singleton or a non-compact space. Due to the following decomposition
    \begin{equation}
        J^{-1}(-\infty,\gamma]=\bigcup_{\tilde{K}\in J_0^{-1}(-\infty,\gamma]} T^{-1}(\tilde{K}),
    \end{equation}
 the sublevel set of $J$ is noncompact in general,  which is the main difficulty in establishing convergence for gradient flows in PO.  In contrast, the coercivity of $J_\lambda$ indicates the compactness of any sublevel set of  $J_\lambda$, which forces the convergence of gradient flows in PO.

From the expressions in Problems $\mathbf{P}_2$ and $\mathbf{P}_3$, the differentiability and analyticity of $J(L,K)$ and $J_\lambda$ at each $(L,K)\in\mathcal{F}$ can be verified straightforwardly. Then we provide the closed form of the policy gradients. Viewing the product space $\mathcal{L}\times\mathcal{K}$ as a smooth manifold, the feasible set $\mathcal{F}$ is an open subset of  $\mathcal{L}\times\mathcal{K}$. Based on \cite[Proposition 3.9]{John2012}, we know that $\mathcal{F}$ is an open submanifold of $\mathcal{L}\times\mathcal{K}$. Given $(L,K)\in\mathcal{F}$,  the following  linear isometry for the \emph{tangent space} holds
\begin{equation}
    \mathcal{T}_{(L,K)}\mathcal{F}\cong \mathcal{L}\times\mathcal{K}\cong \mathbb{R}^{n_l}\times\mathbb{R}^{n_k}.
\end{equation}

For mapping $T:\mathcal{F}\to\mathcal{M}$, we investigate its \emph{tangent map} (or \emph{differential}) as follows when $\mathcal{M}$ admits a local manifold structure\footnote{ Without regularity of $\mathcal{M}$, we can only consider the tangent cone  at $p\in\mathcal{M}$ instead of the tangent space.}. Detailed definitions of the tangent map and other concepts in manifolds can be found in  \cite[Chapter 3]{John2012}.
\begin{lemma}\label{caldt}
Under Assumption \ref{ass:den}, given $\tilde{K}\in \mathbb{R}^{mN\times nN}\setminus S_0$, let $(L,K)\in T^{-1}(\tilde{K})$, then the tangent map $DT_{(L,K)}$ is given by
         \begin{equation}
         DT_{(L,K)}(\Delta L,\Delta K)=(I-L)^{-1}(\Delta K+\Delta L\tilde{K}).
    \end{equation}
  Furthermore, we have
    \begin{equation}\label{dtt}
\begin{aligned}
        \ker(DT_{(L,K)})&=\bigg\{(\Delta L,\Delta K)\in\mathcal{L}\times\mathcal{K}:\Delta K+\Delta L\tilde{K}=0\bigg\},\\
      \rank(DT_{(L,K)})&=\rank \bigg[\hat{K},(\tilde{K}^\top\otimes I_{mN})\hat{L}\bigg].
\end{aligned}
    \end{equation}
\end{lemma}
\begin{proof}  See Appendix \ref{appdis1.5}.
\end{proof}

Endowing the tangent space $\mathcal{T}_{(L,K)}(\mathcal{L}\times\mathcal{K})$
with the Frobenius product $\langle\cdot,\cdot\rangle$ makes
$\mathcal{L}\times\mathcal{K}$ and $\mathcal{F}$ \emph{Riemannian manifolds}.
Under Assumption~\ref{ass:den}, the same metric renders
$\mathcal{M}\cap\mathcal{B}(\tilde{K},\delta)$ a Riemannian manifold
for $\tilde{K}\in\mathbb{R}^{mN\times nN}\setminus S_0$ and 
sufficiently small $\delta>0$.
 Based on the defined metrics, we can introduce the \emph{conjugate map} $DT^\ast_{(L,K)}:\mathcal{T}_{\tilde{K}}\mathcal{M}\to \mathcal{T}_{(L,K)}\mathcal{F}$ of  $DT_{(L,K)}:\mathcal{T}_{(L,K)}\mathcal{F}\to \mathcal{T}_{\tilde{K}}\mathcal{M}$ by
\begin{equation}
\langle (\Delta L,\Delta K),DT^\ast_{(L,K)}W\rangle=    \langle DT_{(L,K)}(\Delta L,\Delta K),W\rangle
\end{equation}
for any $W\in \mathcal{T}_{\tilde{K}}\mathcal{M}$, and $(\Delta L,\Delta K)\in\mathcal{T}_{(L,K)}\mathcal{F}$. 

Since $\mathcal{L}\times\mathcal{K}$ is a submanifold of 
$\mathbb{R}^{mN\times mN}\times\mathbb{R}^{mN\times nN}$,  according to \cite[Proposition 3.53]{boumal2020}, the \emph{Riemannian gradient} $\nabla_{\mathcal{F}} J(L,K)$ just equals the projection of the Euclidean gradient, i.e.,
\begin{equation}\label{gradf}
\nabla_{\mathcal{F}} J(L,K)=
\bigg[\bigg(\Pi_{\mathcal{L}} (\frac{ \partial J}{\partial L})\bigg)^\top,
\bigg(\Pi_{\mathcal{K}} (\frac{ \partial J}{ \partial K})\bigg)^\top\bigg]^\top(L,K),
\end{equation}
where $\Pi_{\mathcal{L}}: \mathbb{R}^{mN\times mN}\to \mathcal{L}$ and $\Pi_{\mathcal{K}}: \mathbb{R}^{mN\times nN}\to \mathcal{K}$ are two natural projection maps\footnote{If $S_\mathcal{L}=\Struct(\mathcal{L})$ and $S_\mathcal{K}=\Struct(\mathcal{K})$, then $\Pi_{\mathcal{L}}(L)=L\odot S_\mathcal{L}$, $\Pi_{\mathcal{K}}(K)=K\odot S_\mathcal{K}$, where '$\odot$' is the Hadamard product.}.  Here $(L,K)$ is called a \emph{stationary point} of $J$ if $\nabla_{\mathcal{F}} J(L,K)=0$. For convenience, we denote two components on the right hand side of (\ref{gradf}) as $\nabla_{\mathcal{L}}J$ and $\nabla_{\mathcal{K}}J$, respectively.
Now we give the expression of gradients and Hessians of cost functions $J$ and $J_\lambda$. 
\begin{Proposition}[Gradients and Hessians of $J$ and $J_\lambda$]\label{pro:grad}
    For $(L,K)\in \mathcal{F}$, we have

    1) The Riemannian gradients of $J$ and $J_\lambda$ at $(L,K)$ are
    \begin{equation}\label{eq:rgra}
    \begin{aligned}
        \nabla_{\mathcal{F}} J(L,K)&=\begin{bmatrix}
\Pi_{\mathcal{L}}\bigg((I-L)^{-\top}\nabla J_0(\tilde{K})\tilde{K}^\top\bigg)\\
\Pi_{\mathcal{K}}\bigg((I-L)^{-\top}\nabla J_0(\tilde{K})\bigg)
\end{bmatrix},\\
        \nabla_{\mathcal{F}} J_\lambda(L,K)&=\nabla_{\mathcal{F}} J(L,K)+\lambda\begin{bmatrix}
L\\
\mathbf{0}
\end{bmatrix}, 
    \end{aligned}
    \end{equation}
 respectively, where $\tilde{K}=T(L,K)$, $\nabla J_0(\tilde{K})=2[R\tilde{K}-B^\top X(A-B\tilde{K})]Y$, and $Y$ is the unique solution to 
\begin{equation}\label{eq:lya}
    (A-B\tilde{K}) Y(A-B\tilde{K})^\top+\Sigma=Y.
\end{equation}

2) Under Assumptions \ref{ass:sol} and \ref{ass:den}, let 
 $\tilde{K}_*\in \mathbb{R}^{mN\times nN}\setminus S_0$ be the optimal feedback gain of $\mathbf{P}_0$, then each
$(L_*,K_*)\in T^{-1}(\tilde{K}_*)$ attains the global minimum of $J$ over $\mathcal{F}$, and near $(L_\ast,K_\ast)$, the Riemannian gradient of $J$ is 
\begin{equation}\label{eq:cha}
   \nabla_{\mathcal{F}} J(L,K)=DT^\ast_{(L,K)}\circ \nabla J_0(\tilde{K}).
\end{equation}
Furthermore, the Hessians of $J$ and $J_\lambda$ at  $(L_\ast,K_\ast)$ are 
\begin{equation}\label{eq:leb}
\begin{aligned}
    \nabla_{\mathcal{F}}^2 J(L_\ast,K_\ast)&= DT_{(L_\ast,K_\ast)}^\ast\circ \nabla^2_{\tilde{K}} J_0(\tilde{K}_\ast)\circ DT_{(L_\ast,K_\ast)},\\
    \nabla_{\mathcal{F}}^2 J_\lambda(L_\ast,K_\ast)
    &= \nabla_{\mathcal{F}}^2 J(L_\ast,K_\ast)+\lambda  \mathcal{Q}_{\mathcal{L}},
    \end{aligned}
\end{equation}
respectively, where $\mathcal{Q}_{\mathcal{L}}$ is a quadratic form on $\mathcal{L}\times\mathcal{K}$ defined by
$\mathcal{Q}_{\mathcal{L}}(\Delta E_1,\Delta E_2)=\trace[(\Delta L_1)^\top \Delta L_2]$  with $\Delta L_i=\Pi_{\mathcal{L}}(\Delta E_i)\in\mathcal{L}$ for $i=1,2$.
\end{Proposition}
\begin{proof}  See Appendix \ref{appdis1.5}.
\end{proof}



\subsection{Projected gradient flows for $J$ and $J_\lambda$: divergence versus convergence}
From (\ref{eq:rgra}), we recall that the Riemannian gradients $\nabla_{\mathcal{F}}J$ and $\nabla_{\mathcal{F}}J_\lambda$ are projections of Euclidean gradients.  In this subsection, we will investigate the convergence behaviors of projected gradient flows (PGFs) for $J$ and $J_\lambda$, and establish theoretical guarantee. 

Note that the PGFs for $J$ and $J_\lambda$ are expressed as
\begin{equation}\label{eq:gfd}
[\dot{L}^\top,\dot{K}^\top]^\top = -\nabla_{\mathcal{F}}J(L,K), 
\end{equation}
\begin{equation}\label{eq:gfr}
[\dot{L}^\top,\dot{K}^\top]^\top = -\nabla_{\mathcal{F}}J_\lambda(L,K),
\end{equation}
respectively, where the expressions of $\nabla_{\mathcal{F}}J$ and $\nabla_{\mathcal{F}}J_\lambda$ are provided in Proposition \ref{pro:grad}. Accordingly, the  PGF for $J_d$ is 
\begin{equation}\label{eq:gfrd}
\dot{K} = -\nabla_{\mathcal{K}}J_d(K),
\end{equation}
where $\nabla_{\mathcal{K}}J_d=\Pi_\mathcal{K}(\nabla J_d)$.

By using the orthonormal bases of $\mathcal{K}$ and $\mathcal{L}$, we can obtain an equivalent representation $\tilde{J}:\Omega\to\mathbb{R}$ of $J(L,K)$, defined by $\tilde{J}(x,y)=J(\vvec^{-1}(\hat{L}x),\vvec^{-1}(\hat{K}y))$, where  $\Omega\triangleq\{(x,y)\in\mathbb{R}^{n_l}\times\mathbb{R}^{n_k}: (\vvec^{-1}(\hat{L}x),\vvec^{-1}(\hat{K}y)\in\mathcal{F}\}$. Then $\nabla_{\mathcal{F}} J(L,K)$ can be identified with the usual gradient on $\mathbb{R}^{n_l}\times\mathbb{R}^{n_k}$ without constraints. Therefore, results on gradient flows in Euclidean spaces  (see, e.g., \cite{2005gra,stab2006,2025MP}) can be applied directly to our setting. 




First, we show that given any initial condition $(L_0,K_0)\in\mathcal{F}$, the long time existence and uniqueness of solutions hold for both Equations (\ref{eq:gfd}) and (\ref{eq:gfr}).
\begin{lemma}[Long-time existence of solutions to PGFs] \label{long}
    Given any initial condition $(L_0,K_0)\in\mathcal{F}$, there exists a unique solution $(L(t),K(t))\in C^\infty([0,+\infty),\mathcal{F})$ to (\ref{eq:gfd}) with $L(0)=L_0$ and $K(0)=K_0$. Moreover, the same conclusion holds for (\ref{eq:gfr}).
\end{lemma}
\begin{proof}
See Appendix \ref{appdis2}.
\end{proof}

\begin{figure*}[!htbp]
\centering
\includegraphics[width=\linewidth]{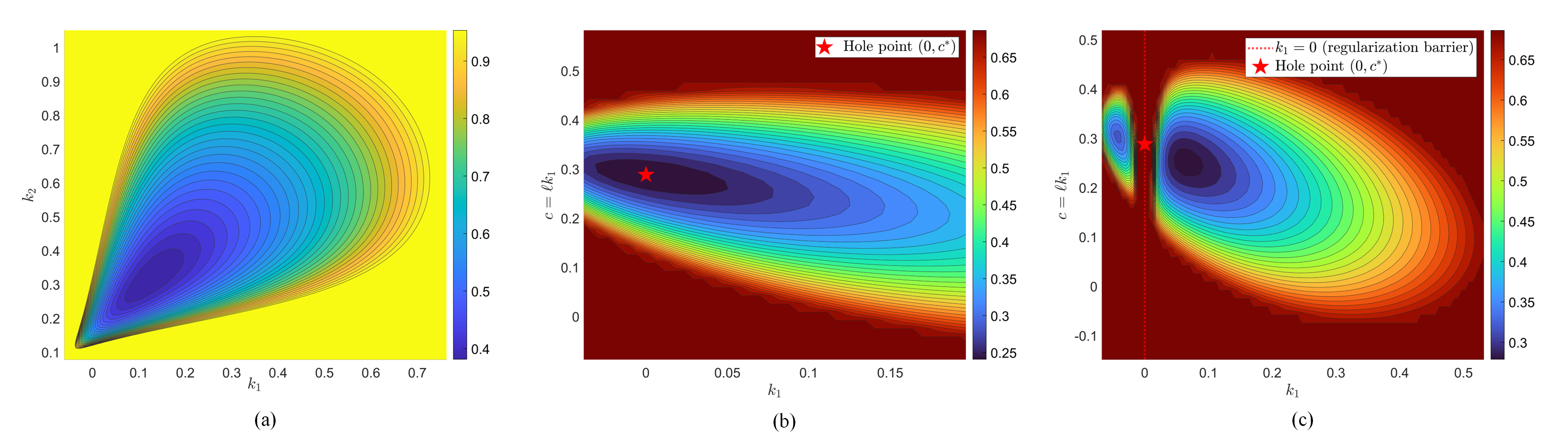}
\caption{Function landscapes of $J_d$, $J$, and $J_\lambda$ with $\lambda=0.001$. (a) $k_1$-$k_2$ plot, where the minimal cost $J_{d,\ast}=0.3813$ is achieved at $K_*=\diag(0.1239, 0.3147)$. (b) $k_1$-$c$ plot  on the slice $k_2
 =0.6601$, where the infimum of cost $J_\ast=0.2404$ is not attainable since it corresponds to a hole point $(k_1^\ast,k_2^\ast,c^\ast)=(0, 0.6601,0.2881)$. (c) $k_1$-$c$ plot on the same slice $k_2
 =0.6601$, where the regularization term induces a barrier along $k_1=0$. The minimal cost $J_{\lambda,\ast}=0.2745$ is achieved at $(k_1^\ast,k_2^\ast,c^\ast)=( 0.056, 0.561,0.199)$, which lies off the plotted slice. }  
\label{fig:landsc}
\end{figure*}
\begin{figure*}[htbp]
\centering
\includegraphics[width=\linewidth]{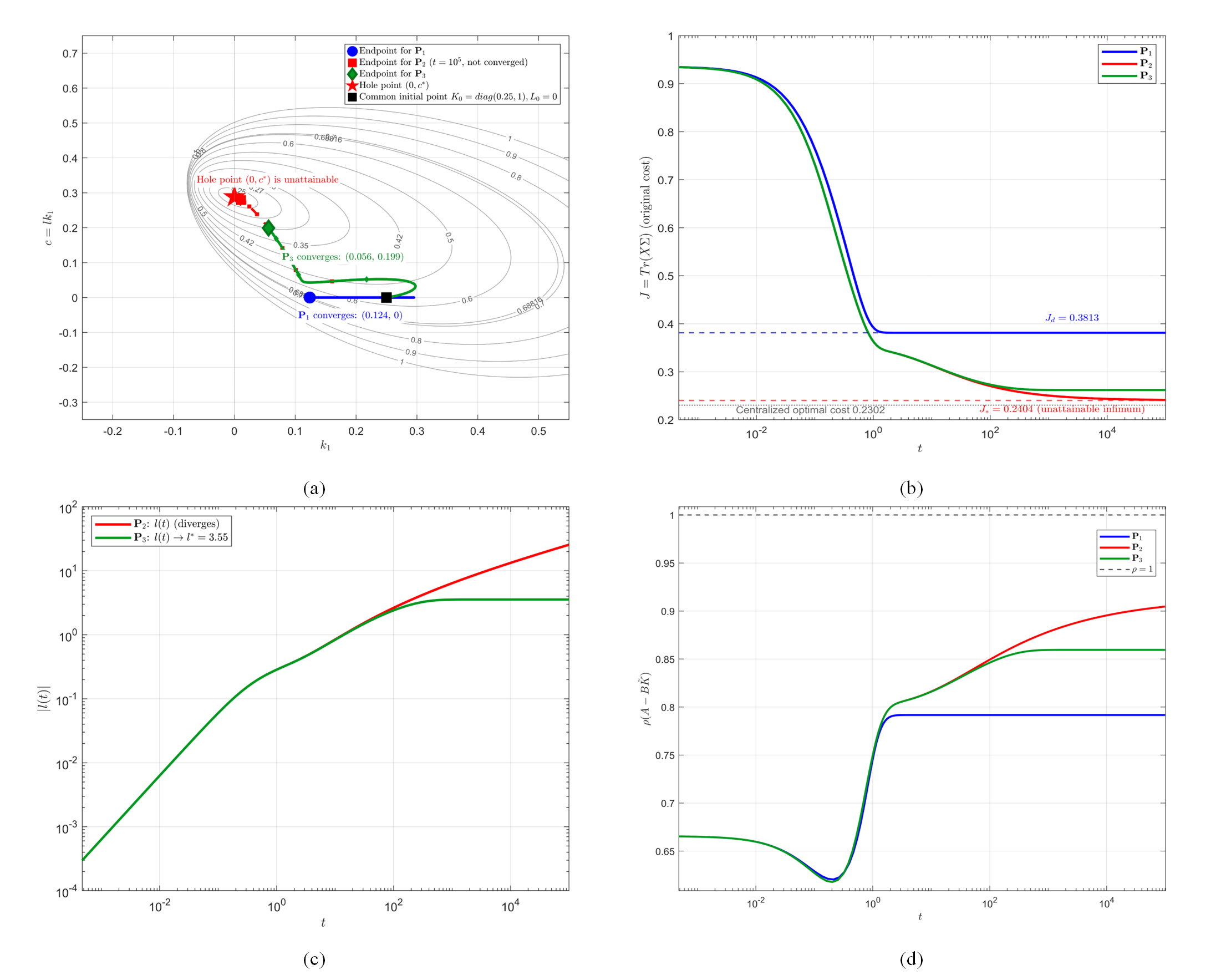}
\caption{Convergence behaviors for PGFs of $J_d$, $J$, and $J_\lambda$ with the same initial point $K_0$ or $(0,K_0)$ and $\lambda=0.001$. (a) The evolution of trajectories on the $k_1$-$c$ plane along three PGFs. Note that the hole point corresponding to $(0,c^\ast)$ is not achievable along the PGF of $J$. (b) The evolution of the original cost. (c) The evolution of the norm $|\ell(t)|$. It is noteworthy that $|\ell(t)|$ is unbounded along the PGF of $J$. (d) The evolution of the spectral radius for the closed-loop system.}
\label{fig:grade}
\end{figure*}
One might intuitively expect the solutions of (\ref{eq:gfd}) to converge for arbitrary initial conditions in $\mathcal{F}$. Now we show via a counterexample that it is not true in general. Furthermore, we reveal how the optimization landscapes of Problems  
$\mathbf{P}_1$, $\mathbf{P}_2$ and $\mathbf{P}_3$
   lead to distinct convergence behaviors of PGFs. \\
\noindent\textbf{Example 2:} [Counterexample to the convergence of solutions to (\ref{eq:gfd})]\label{ex:2node}
Let $A = \diag(a_1, 1.05)$, $B = R = I_2$, $m=n=1$, and $N=2$, where $a_1$ is a constant to be determined.
Take 
$$Q = \Bigg[\begin{matrix} 0.3 & \tfrac{199\sqrt{30}}{2000} \\ \tfrac{199\sqrt{30}}{2000} & 1.0
\end{matrix}\Bigg],\ \Sigma = \Bigg[\begin{matrix} 1 & -\frac{9\sqrt{2}}{20} \\ -\frac{9\sqrt{2}}{20} & 0.5
\end{matrix}\Bigg],\ \mathcal K = \Bigg\{\diag(k_1, k_2): k_1,k_2\in\mathbb{R} \Bigg\}, \mathcal{L}=\Bigg\{\Bigg[\begin{matrix} 0 & 0 \\  \ell & 0
\end{matrix}\Bigg]:\ell\in\mathbb{R}\Bigg\}.$$
 Then we can verify that $\mathcal{M}$ is
 \begin{align}
\mathcal M 
 = &\Bigg\{ \Bigg[\begin{matrix} k_1 & 0 \\ c & k_2
\end{matrix}\Bigg]: k_1,k_2,c\in\mathbb{R}\Bigg\}
\Bigm\backslash \Bigg\{ \Bigg[\begin{matrix} 0 & 0 \\ c & k_2
\end{matrix}\Bigg]: k_2\in\mathbb{R},c\in\mathbb{R}\setminus\{0\}\Bigg\},
\end{align}
and the closure of $\mathcal{M}$ is
$$\bar{\mathcal M} = \Bigg\{ \Bigg[\begin{matrix} k_1 & 0 \\ c & k_2
\end{matrix}\Bigg]: k_1,k_2,c\in\mathbb{R}\Bigg\}.$$ 
Note that $J(L,K)$ can be viewed as a function of $(k_1,k_2,c)\in \mathbb{R}^3$. Specifically, it is quadratic with respect to  $c$. After minimizing with respect to $c$, we get a function of $(k_1,k_2)$ with parameter $a_1$. If we fix $k_1 = 0$, the
profile $c^\ast(k_2)$ is available in closed form, and 
minimization over $k_2$  yields the 
minimizer $(0, k_2^\ast, c^\ast)$ over the set $\{(0,k_2,c):k_2,c\in\mathbb{R}\}$. Then solving $\frac{\partial J}{\partial k_1}\Big|_{(0,\ k_2^\ast,\ c^\ast)}
  \;=\; 0$ yields a choice of $a_1\approx 0.9156$. 

Figure \ref{fig:landsc} presents landscapes of cost functions $J_d$, $J$, and $J_\lambda$. Based on the contour plot map of $J$ in Figure \ref{fig:landsc}(b), we can obtain that $$\tilde{K}_\ast \approx \Bigg[\begin{matrix} 0 & 0 \\ 0.2881 & 0.6601
\end{matrix}\Bigg]\in\bar{\mathcal{M}}\setminus\mathcal{M}$$ corresponds to the infimum of $J$, and is in fact a \emph{hole point} \footnote{Here $\tilde{K}_*$ is called a hole point if there exists a sequence $\{(L_j,K_j)\}\subseteq\mathcal{F}$  such that $\lim_{j\to\infty}T(L_j,K_j)=\tilde{K}_*$, while $\tilde{K}_*\notin\mathcal{M}$.
}.  As shown in Figure \ref{fig:landsc}(c), the partially regularized term drives the cost up sharply near the regularization barrier  $k_1=0$. 

Figure \ref{fig:grade} illustrates different convergence behaviors of PGFs. From Figure \ref{fig:grade}(a), we observe different trajectories of solutions to PGFs.  Figure \ref{fig:grade}(b) shows an obvious improvement of performance of $\mathbf{P}_2$ and $\mathbf{P}_3$ compared with $\mathbf{P}_1$. Figure \ref{fig:grade}(c) illustrates the divergence  of $\ell(t)$ along the PGF of $J$, while the convergence holds for that of $J_\lambda$. Figure \ref{fig:grade}(d) indicates that along PGFs, the performance improvement comes at the cost of  a weaker stability of the closed-loop system. Compared with $\mathbf{P}_2$,  $\mathbf{P}_3$ achieves a better balance between stability and performance.                      \hfill $\square$

The above example suggests that no convergence guarantee can be
expected for the PGF of $J$ without further assumptions.
Fortunately, for the PGF  of $J_\lambda$, the convergence holds for free due to the coercivity of $J_\lambda$.
\begin{theorem}\label{thm:reg}
     For any initial condition $(L_0,K_0)\in\mathcal{F}$, there always exists a stationary point $(L_\ast,K_\ast)\in\mathcal{F}$ of $J_\lambda$ such that the solution $(L(t),K(t))$ to (\ref{eq:gfr}) 
     converges to $(L_\ast,K_\ast)$. Furthermore, the following asymptotic convergence rate holds 
     \begin{equation}
     \begin{aligned}
           J_\lambda(L(t),K(t))-J_\lambda(L_\ast,K_\ast)&=\mathcal{O}(t^{-\frac{1}{2\mu-1}}), \
         \|(L(t),K(t))-(L_\ast,K_\ast)\|=\mathcal{O}(t^{-\frac{1-\mu}{2\mu-1}}), \text{if}\ \mu>\frac{1}{2};\\
            J_\lambda(L(t),K(t))-J_\lambda(L_\ast,K_\ast)&=\mathcal{O}(\exp(-c^2t)),\
         \|(L(t),K(t))-(L_\ast,K_\ast)\|=\mathcal{O}(\exp(-c^2 t/2)), \text{if}\ \mu=\frac{1}{2},
     \end{aligned} 
     \end{equation}
     where $\mu>0$ is the Łojasiewicz exponent, and $c>0$ is the constant in Łojasiewicz inequality of $J_\lambda$ at $(L_\ast, K_\ast)$.
 \end{theorem}  

\begin{proof}
 See Appendix \ref{appdis4}.
\end{proof}

    

 \section{Globally optimal policy synthesis}\label{optim}
 
In this section, we first bridge Problems $\mathbf{P}_1$ and $\mathbf{P}_2$ via the partially regularized Problem $\mathbf{P}_3$. Then we investigate how to design optimal policies by PO when the centralized optimal policy is representable under action dependencies.  Finally, the improvement of performance under general ADGs will be briefly discussed. 

Let $\tilde{K_\ast}$ be the centralized optimal feedback gain of Problem $\mathbf{P}_0$. 
To guarantee $\tilde{K_\ast}\in\mathcal{M}\cap\mathcal{S}$, we propose the following assumption.  
\begin{assumption}\label{ass5}
Suppose that Assumptions \ref{ass:st} and \ref{ass:den}  hold, and $\tilde{K_\ast}\in\mathbb{R}^{mN\times nN}\setminus S_0$ is the centralized optimal feedback gain corresponding to the unique solution of ARE (\ref{eq:are}).   
\end{assumption}

Based on \cite[Theorem 2.4.2]{lewis2012optimal}, Assumption \ref{ass5} indicates that the centralized optimal policy $\tilde{K_\ast}$ stabilizes the original system. Together with Theorem \ref{thm:c1},  we have $\tilde{K_\ast}\in\mathcal{M}\cap\mathcal{S}$ and $T^{-1}(\tilde{K_\ast})\subseteq\mathcal{F}$. Then the desired equilibrium set is defined as
$ \mathcal{E}_\ast\triangleq T^{-1}(\tilde{K_\ast})$,
which actually corresponds to the set of global minimizers. 
\subsection{Connection between Problems $\mathbf{P}_1$ and $\mathbf{P}_2$ via partial regularization}
Let  $J_{d,\ast}$, $J_\ast$, and $J_{\lambda,\ast}$ be the minimal costs of Problems $\mathbf{P}_1$, $\mathbf{P}_2$ and $\mathbf{P}_3$, respectively. Furthermore, the set of distributed controllers achieving the minimal cost $J_{d,\ast}$ is denoted by $ 
\hat{\mathcal{E}}_\ast\triangleq\{\hat{K}_\ast\in\mathcal{K}:J_d(\hat{K}_\ast)=J_{d,\ast}\}$. For convenience, let $(L_\lambda,K_\lambda)\in\mathcal{F}$ be a global minimizer of $J_\lambda$ (may be non-unique), i.e., $J_{\lambda}(L_\lambda,K_\lambda)=J_{\lambda,\ast}$. 

Under Assumption \ref{ass5}, 
we denote
    \begin{align}\label{delmi}
        \bar{L}_\ast
        \triangleq\argmin\bigg\{\|L\|^2: T(L,K)=\tilde{K}_\ast, (L,K)\in\mathcal{L}\times\mathcal{K}\bigg\},\ 
        \bar{K}_\ast\triangleq(I-\bar{L}_\ast)\tilde{K}_\ast.
    \end{align}
        To illustrate the definition of $\bar{L}_*$ in (\ref{delmi}),  we identify the fiber 
    $T^{-1}(\tilde{K}_\ast)=\{(L,K)\in\mathcal{L}\times\mathcal{K}:K=(I-L)\tilde{K}_\ast\}$
   with set $\mathcal{C}_{\mathcal{L}}\triangleq\{L\in\mathcal{L}: (I-L)\tilde{K}_\ast\in\mathcal{K}\}$ via the map $L\mapsto\bigl(L,(I-L)\tilde K_\ast\bigr)$. Note that $\mathcal{C}_{\mathcal{L}}$ is nonempty, closed and convex.  According to \cite[Theorem 2.2.11]{nes2018}, it has a unique element of minimal norm, denoted by $\bar{L}_\ast\in\mathcal{L}$.
Correspondingly, $\bar{K}_\ast=(I-\bar{L}_\ast)\tilde{K}_\ast$. Thus $\bar{L}_\ast$ in (\ref{delmi}) is well-defined.

 \begin{theorem}\label{conne}\textit{(Bridging  Problems $\mathbf{P}_1$ and $\mathbf{P}_2$ via partial regularization)}
  Under Assumption \ref{ass:st}, the following convergence properties hold: 
    \begin{align}
        \lim_{\lambda\to 0}J_{\lambda,\ast}=J_\ast,\ \lim_{\lambda\to\infty}J_{\lambda,\ast}=J_{d,\ast},\
        \lim_{\lambda\to \infty}L_\lambda=0,\ \lim_{\lambda\to \infty}\dist(K_\lambda,\hat{\mathcal{E}}_\ast)=0.
    \end{align}
    Furthermore, under Assumption \ref{ass5},  the following properties hold as well:
    \begin{equation}\label{con}
        \lim_{\lambda\to 0}L_\lambda=\bar{L}_\ast,\ \lim_{\lambda\to 0}K_\lambda=\bar{K}_\ast.
    \end{equation}
\end{theorem}
\begin{proof}
    See Appendix \ref{bridge}.
\end{proof}

Theorem \ref{conne} is a statement about the connection between $\mathbf{P}_1$ and $\mathbf{P}_2$, not about specific optimization algorithms. It indicates that under Assumption \ref{ass5}, $(L_\lambda,K_\lambda)$ always converges to the unique point $(\bar{L}_\ast,\bar{K}_\ast)$ with minimal $L$-component norm on the fiber $T^{-1}(\tilde{K}_*)$ as $\lambda\to0$. Interestingly, as $\lambda\to\infty$, $L_\lambda\to 0$, while $K_\lambda$ converges to the set of minimizers of $\mathbf{P}_1$. This matches the intuition behind partial regularization: a larger $\lambda$ penalizes the action-dependent term more heavily, forcing the norm of $L$ to shrink. 

Based on (\ref{pr11}) and (\ref{pr2}), we provide a concise description of some convergence rates as follows.
\begin{corollary}\label{cor5}
 Under Assumption \ref{ass:st}, the following convergence rates hold
\begin{align}
 J_{\lambda,\ast}-J_{\ast}=\mathcal{O}(\lambda)\ \text{as}\ \lambda\to 0,\   \|L_{\lambda}\|=\mathcal{O}(\lambda^{-\frac{1}{2}})\ \text{as}\ \lambda\to \infty.
\end{align}
\end{corollary}
\subsection{Convergence analysis of the PGFs for costs $J$ and $J_\lambda$}
Recall that the PGF of cost $J$ may fail to converge in general as illustrated in Example 2. In this subsection, we show that Assumption \ref{ass5} plays a crucial role in ensuring the convergence of the PGF near optimal policies for Problem $\mathbf{P}_2$. Analogous result also holds for  the PGF of cost $J_\lambda$.

Note that under Assumption \ref{ass5}, $\mathcal{E}_\ast$ is a submanifold (affine subspace) of $\mathcal{L}\times\mathcal{K}$, i.e., the equilibria are usually non-isolated. 
Viewing $T:\mathcal{F}\to\mathcal{M}\cap\mathcal{S}\subseteq\mathbb{R}^{mN\times nN}$ as a parameterization mapping, previous discussion in Remark \ref{re:ass} shows that $T$ is not injective and the dimension of fiber $T^{-1}(\tilde{K})$ may vary as $\tilde{K}$ varies. This leads to main difficulties in understanding global landscapes of $J$.

Nevertheless, we prove that, under Assumption \ref{ass5}, $T$  has constant rank in a neighborhood of  
$\tilde{K}_*$, so that the fibers $T^{-1}(\tilde{K})$ are mutually isometric for  $\tilde{K}$ in that neighborhood. 

\begin{Proposition}\label{sub}
Under Assumption \ref{ass5}, $T$ is a submersion in a neighborhood of $ \mathcal{E}_\ast$ for generic $\tilde{K}_*$. 
\end{Proposition}
\begin{proof}
    See Appendix \ref{psub}.
\end{proof}

Then we turn to investigating the finer optimization property of $J$ at $(L_\ast,K_\ast)\in\mathcal{E}_\ast$. Before that, we recall the definition of \emph{Morse-Bott property} (MB) for a smooth function. 
\begin{definition}\cite[Morse--Bott property at a local minimum]{2025MP}
\label{MBde}
Let $\mathcal{X}$ be a smooth manifold, and $f\in C^\infty(\mathcal{X})$. Let $p_\ast$ be a local minimum of $f$ with associated set $\mathcal{O}$ defined by $$\mathcal{O}=\{p\in\mathcal{X}: f(p)=f(p_\ast)\}.$$ We say $f$ satisfies the Morse--Bott property at $p_\ast$ if
$\mathcal{O}$ is a  $C^1$ submanifold around $p_\ast$ and $\ker \nabla^2 f(p_\ast) = \mathcal{T}_{p_\ast}\mathcal{O}$.
If $\langle v, \nabla^2 f(p_\ast) [v] \rangle \geq \mu \|v\|^2$ also holds
for some $\mu > 0$ and all $v \in \mathcal{N}_{p_\ast} \mathcal{O}$, where $\mathcal{N}_{p_\ast} \mathcal{O}$ is the orthogonal complement of $\mathcal{T}_{p_\ast}\mathcal{O}$ in $\mathcal{T}_{p_\ast}\mathcal{X}$, then we say that $f$ satisfies the $\mu$-MB property at $p_\ast$.
\end{definition}

Intuitively,  the MB property of $f$ at $p_*$ indicates that the
Hessian of $f$ at $p_*$ is positive definite in the normal directions. Now we reveal the MB property of $J$ at $(L_\ast,K_\ast)\in\mathcal{E}_\ast$, and then obtain \emph{Polyak–Łojasiewicz} (PL) inequality and \emph{quadratic growth} (QG) property of $J$ in a neighborhood of $(L_\ast,K_\ast)$.
\begin{Proposition}\label{lem:MBP}
     Under Assumption \ref{ass5}, suppose $(L_\ast,K_\ast)\in\mathcal{E}_\ast$, denote the smallest positive eigenvalue of $\nabla^2_{\mathcal{F}}J(L_\ast,K_\ast)$ as $\lambda_\ast\triangleq\lambda_\ast(L_\ast,K_\ast)$, then the following statements are true:
    
     1) $J$ has $\lambda_*$-MB property at  $(L_\ast,K_\ast)$. 

     2) there exist  constants $\mu,\mu^\prime\in(0,\lambda_\ast)$ and neighborhood $\mathcal{U}$ of  $(L_\ast,K_\ast)$ such that 
  \begin{itemize}  
      \item PL inequality:
    \begin{equation}\label{muPL}
        J(L,K)-J_\ast \leq \frac{1}{2\mu}\|\nabla_{\mathcal{F}}J(L,K)\|^2
    \end{equation}
      \item QG property: 
        \begin{equation}\label{muQG}
        J(L,K)-J_\ast \geq \frac{\mu^\prime}{2}\dist[(L,K),\mathcal{E}_\ast]^2
    \end{equation}
    hold for $(L,K)\in \mathcal{U}$, and $J_\ast=J(L_\ast,K_\ast)$.
  \end{itemize}
\end{Proposition}
\begin{proof}
 See Appendix \ref{appdis5}.
\end{proof}

Building on Proposition \ref{lem:MBP}, we establish the asymptotic stability of $\mathcal{E}_*$ as follows.


    

\begin{theorem} \label{thm:c3} \textit{(Local exponential convergence to global optimality of $J$)}
Under Assumption \ref{ass5}, for any $(L_\ast,K_\ast)\in\mathcal{E}_\ast$, there exists a neighborhood $\mathcal{U}_0$ of $(L_\ast,K_\ast)$ such that for initial condition $(L_0,K_0)\in \mathcal{U}_0$, the PGF of $J$ satisfies

1) $(L(t),K(t))$  converges to some $(L_\infty,K_\infty)\in\mathcal{E}_\ast$  as $t\to\infty$;

2) both $\dist[(L(t),K(t)),\mathcal{E}_\ast]$ and  $J(L(t),K(t))-J_\ast$ converge to $0$ exponentially  as $t\to\infty$. 
\end{theorem}
 \begin{proof}
 See Appendix \ref{appdist2}.
 \end{proof}

According to Proposition \ref{pron}, a special case occurs, i.e.,  $\mathcal{E}_\ast$ is a singleton if the following further holds:  
\begin{equation}\label{con1}
n_l+n_k-\rank \bigg[\hat{K},(\tilde{K}_*^\top\otimes I_{mN})\hat{L}\bigg]=0.
\end{equation}
Indeed,  $T$ is a local diffeomorphism near $\mathcal{E}_\ast$ in this setting. 
\begin{corollary}
    Under Assumption \ref{ass5}, if (\ref{con1})  holds, then $\mathcal{E}_\ast$ is a singleton and exponentially asymptotically stable (EAS).  
\end{corollary}





In what follows, we focus on the PGF for Problem $\mathbf{P}_3$. Based on the connection established in Theorem \ref{conne}, we prove that every minimizer $(L_\lambda,K_\lambda)$ of $J_\lambda$   
is EAS. 
\begin{theorem}\label{thm:c6} \textit(Exponential convergence of PGF for $J_\lambda$)
    Under Assumption \ref{ass5}, there exist a neighborhood $\mathcal{U}$ of $(\bar{L}_*,\bar{K}_*)$, which is defined by (\ref{delmi}), and a constant $\delta>0$  such that the following statements hold for $\lambda\in(0,\delta)$:

    1) $(L_\lambda,K_\lambda)\in \mathcal{U}$;
    
    2) $(L_\lambda,K_\lambda)$ is  EAS for the PGF of $J_\lambda$.
\end{theorem}
\begin{proof}
    See Appendix \ref{appthm5}.
\end{proof}

 In fact, due to the partially regularized term, $J_\lambda$ is strongly convex in a neighborhood $\mathcal{U}$ of $(\bar{L}_\ast,\bar{K}_\ast)$, and minimizer $(L_\lambda,K_\lambda)$ must fall into $\mathcal{U}$ for sufficiently small $\lambda$. The above theorem actually indicates that under Assumption \ref{ass5}, the global minimizer of $J_\lambda$ is not only isolated, but also unique when $\lambda$ is sufficiently small. 

\begin{remark}
    Notably, if $\dim(\mathcal{E}_*)>0$, each equilibrium point in $\mathcal{E}_*$ is stable, but not asymptotically stable for the PGF of $J$ since it is not isolated. On the other hand, one may wonder whether $\mathcal{E}_\ast$ is globally asymptotically stable (GAS) under Assumption \ref{ass5}. The answer is negative in general. We can verify that the example in Appendix \ref{appdis1} satisfies Assumptions \ref{ass:st} and \ref{ass:den}. However, the feasible set $\mathcal{F}$ is disconnected, and $\mathcal{E}_\ast$ is always contained in one connected component of $\mathcal{F}$. Therefore,  $\mathcal{E}_\ast$ is not GAS here. Moreover, analogous result also holds for the equilibrium set of the PGF for $J_\lambda$.
\end{remark} 
\subsection{Discussion on the improvement of performance}

Here we discuss the improvement of performance due to ADPs to answer the question proposed in Subsection \ref{prob}.

First, we note that globally optimal policy synthesis in previous subsections enables the improvement of performance. Under Assumption \ref{ass5}, Theorem \ref{thm:c3} implies that the equilibrium set $\mathcal{E}_*$ of (\ref{eq:gfd}) is EAS. Therefore, if initial condition $(L_0, K_0)$ lies in the attraction domain of $\mathcal{E}_*$, then the final suboptimality gap can be reduced to 0, and a sufficient improvement of performance is guaranteed. For the partially regularized case, there exists a bias with an order $\mathcal{O}(\lambda)$ of the optimal cost as $\lambda\to0$ according to Corollary \ref{cor5}.

Then we consider the general case where Assumption \ref{ass5} does not hold and globally optimality cannot be achieved. Specifically, we provide Algorithm \ref{alg:buildpar} to guarantee the improvement of performance compared with the distributed LQR approach.   
\begin{algorithm}
\caption{Algorithm for finding $L$ and $K$ to guarantee the improvement of performance}
\label{alg:buildpar}
\begin{algorithmic}[1]
\STATE{Initialization: $\epsilon>0$, $K_0\in\mathcal{K}\cap\mathcal{S}$}
\STATE{Implement (\ref{eq:gfrd}) until convergence of $K(t)$, i.e., $K(t)\to K_\infty$ as $t\to\infty$}
\STATE{Set $K(0)=K_\infty$ and $L(0)=0$}
\STATE{Implement (\ref{eq:gfd}) until $\|\nabla_{\mathcal{F}}J(L(t),K(t))\|< \epsilon$}
\RETURN $L$ and $K$ 
\end{algorithmic}
\end{algorithm}
Note that for $K_0\in\mathcal{K}\cap\mathcal{S}$, the solution to the PGF of $J_d$ must converge to some stationary point $K_\infty\in\mathcal{K}\cap\mathcal{S}$ of $J_d$. Namely, Line 2 in Algorithm \ref{alg:buildpar} guarantees $K_\infty\in\mathcal{K}\cap\mathcal{S}$ satisfying $ \nabla_{\mathcal{K}}J_d(K_\infty)=0$. Then Lines 3-5 ensure that $(L,K)\in\mathcal{F}$ satisfies $\|\nabla_{\mathcal{F}}J(L,K)\|<\epsilon$ in finite time\footnote{Despite the possible divergence of the solution to (\ref{eq:gfd}), the criterion in Line 4 guarantees finite time stopping. Otherwise, the decrease $J(0,K_\infty)-J(L(t),K(t))=\int^t_0\|\nabla_{\mathcal{F}}J(L(s), K(s))\|^2ds\geq \epsilon^2 t$ indicates that  $\lim\limits_{t\to\infty}[J(0,K_\infty)-J(L(t),K(t))]=\infty$. This contradicts the non-negativity of $J(\cdot,\cdot)$.}. In addition, the following inequality holds naturally

\begin{equation}
    J(L,K)\leq J(0,K_\infty)=J_d(K_\infty)\leq J_d(K_0),
\end{equation} 
 which indicates the improvement of performance compared with the usual distributed LQR. However, how much the performance can be improved, as well as the relationship between the suboptimality gap and ADG, remains an open question and is left for future work. 

Now we revisit Example 2 to illustrate the improvement of performance achieved by Algorithm \ref{alg:buildpar}. We choose $\epsilon=1e-4$, and Line 2 offers $K_\infty=\diag(0.1239, 0.3147)$ and $J_d(K_\infty)=0.3813$. Then by executing Lines 3-5, we obtain that      
$K=\diag(0.0122, 0.6529)$, and $L$ is lower triangular with $L(2,1)= 22.7794 $, and $J(L,K)= 0.2416$. 
Compared with purely executing the PGF of $\mathbf{P}_1$, the improvement of performance is\footnote{Here $J_*$ represents the  optimal cost achieved by solving ARE \eqref{eq:are}.} $\frac{J_d(K_\infty)-J(L,K)}{J_*}=\frac{0.3813-0.2416}{0.2302}\approx 60.67\%$, which is significant.

\section{Numerical simulations}
\label{num}
Consider an MAS consisting of 6 agents with $m=n=2$. Denote $A_0=\begin{bmatrix}
    2 & 1  \\
   1 &  3
\end{bmatrix}$ and 
\begin{equation}
A=\blkdiag\bigg(A_0+ \kappa_{h}\times \text{rand}(2, 2)\bigg)_{6\times 6}+
 \kappa_{c} (\mathbf{1}_6\mathbf{1}^\top_6-I_6)\otimes (\mathbf{1}_2\mathbf{1}^\top_2),
B=\Sigma=I_{12}, R=0.3I_{12},
\end{equation}
where $\kappa_{h}$ is introduced to adjust the heterogeneity of the dynamics,  and rand(2,2) is a random 2-by-2 matrix with elements lying in $[0,1]$, and $\kappa_{c}$ describes the coupling effects among the dynamics of agents. Here we take $\kappa_{h}=1$ and $\kappa_{c}=0.2$. Denote the adjacency matrices of sensing graph $\mathcal{G}_S$ and ADG $\mathcal{G}_A$, and $Q_0$ in (\ref{defgr}) (see Appendix \ref{parae}). Let the state weight matrix be $Q=Q_0\otimes I_2$.

By solving (\ref{eq:are}) numerically, we obtain the centralized optimal controller $\tilde{K}_*$, and the optimal cost of $\textbf{P}_0$ is verified to be $J_*= 92.038$. Choose $K_0\in\mathcal{K}$ as a common initial stabilizing policy.  

\begin{figure*}[htbp]
\centering
\includegraphics[width=0.9\linewidth]{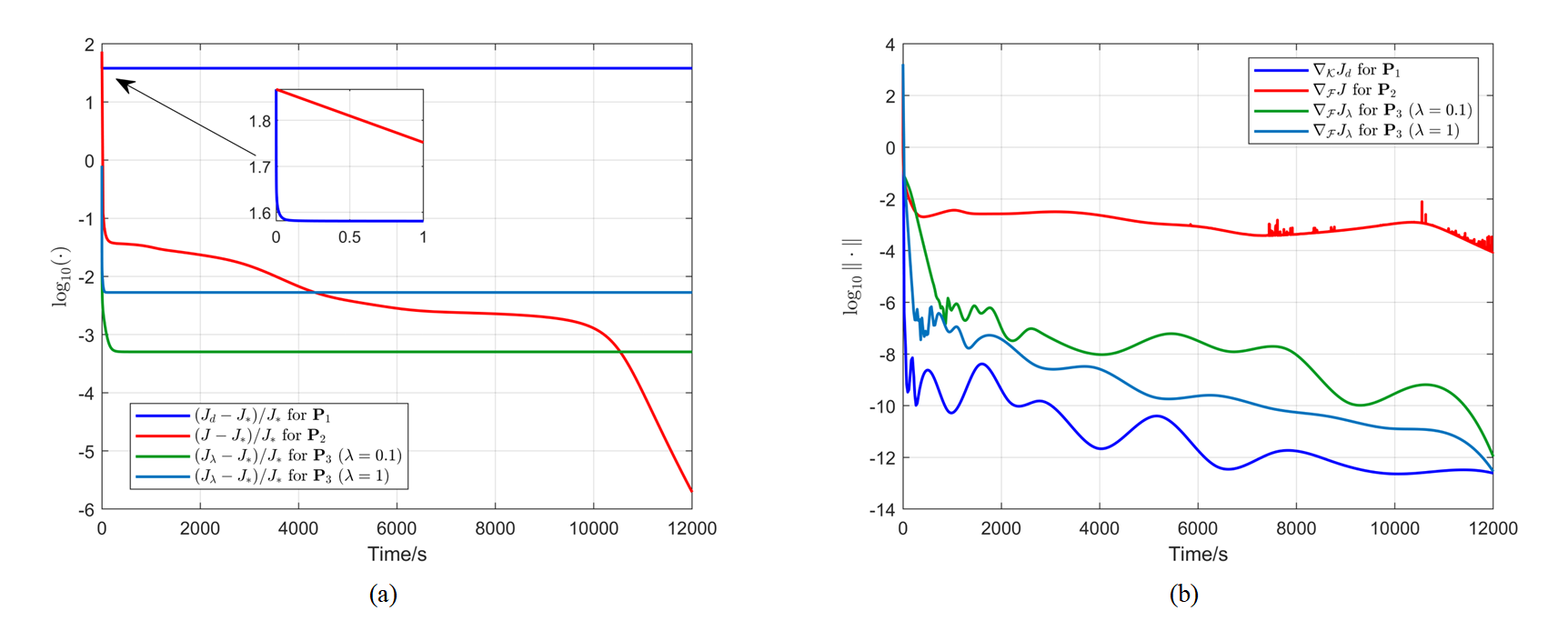}
\caption{The performance evolution of PGFs corresponding to $\mathbf{P}_1$, $\mathbf{P}_2$ and $\mathbf{P}_3$. }
\label{PGF4}
\end{figure*}

\noindent\textbf{Convergence of PGFs and optimality guarantee:} One can verify that Assumption \ref{ass5} holds in this setting. Consider PGFs for $\mathbf{P}_1$ with $K(0)=K_0$, $\mathbf{P}_2$ and $\mathbf{P}_3$ with $L(0)=0$, and $K(0)=K_0$. Figure \ref{PGF4} illustrates the evolution of  $\frac{J-J_*}{J_*}$ and norms of projected gradients. Figure \ref{PGF4}(a) indicates that the solution to the PGF of $\mathbf{P}_2$ finally achieves the global optimality, while that of $\mathbf{P}_1$ only converges to a suboptimal value.  Due to the effects of partial regularization, the inevitable bias of the converged policy for $\mathbf{P}_3$ with respect to the optimal case indeed exists, and a smaller penalty parameter $\lambda$ leads to a smaller bias. This reflects the power of ADPs, and is consistent with Theorems \ref{thm:c3} and \ref{thm:c6}. 
From Figure \ref{PGF4}(b), we observe that the convergence in terms of  $\nabla_{\mathcal{K}} J_d$ is fastest, while that of  $\nabla_{\mathcal{F}} J$ is slowest.  For $\mathbf{P}_3$, a larger penalty parameter $\lambda$ leads to a faster convergence, but also a larger bias with respect to the global optimality. In other words, $\mathbf{P}_3$ offers a trade-off between convergence speed and suboptimality. 

\begin{figure*}[!htbp]
\centering
\includegraphics[width=0.9\linewidth]{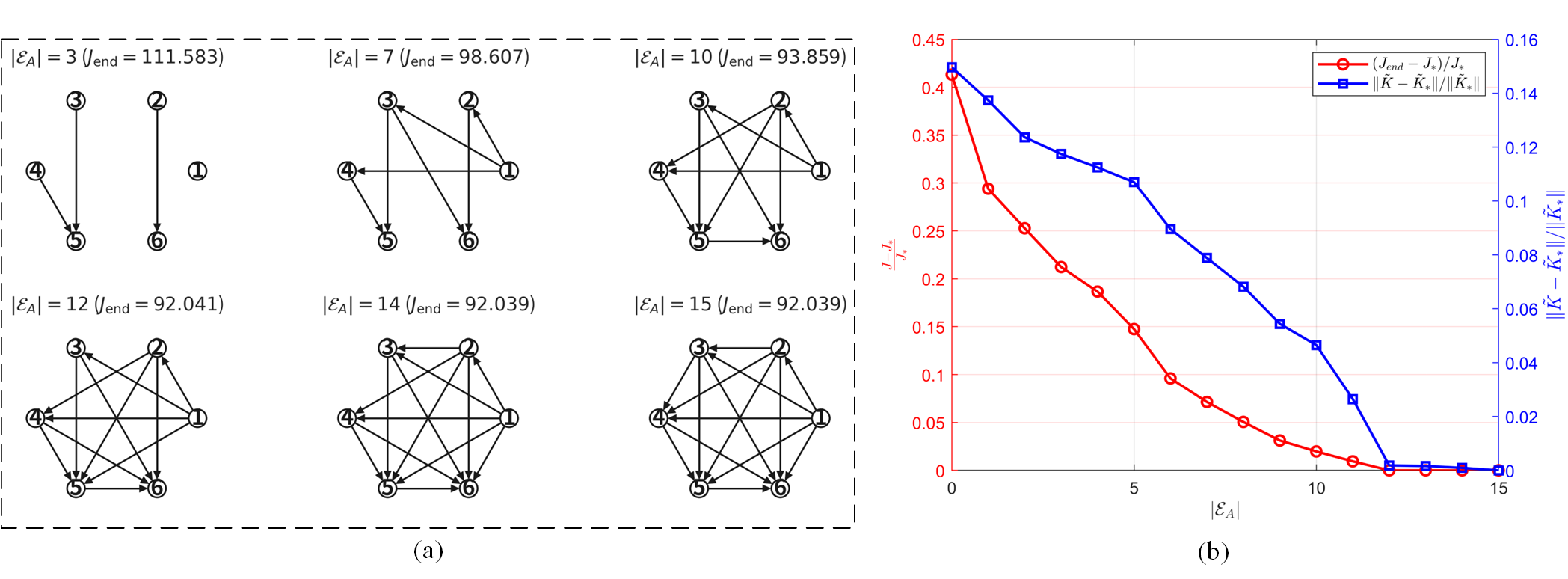}
\caption{Illustration of the effects of ADGs on suboptimality gap. (a) A sequence of ADGs with increasing density. (b) Variations of $\frac{J-J_*}{J_*}$ and $\frac{\|\tilde{K}-\tilde{K}_*\|}{\|\tilde{K}_*\|}$ due to different ADGs. Here $J_{end}$ and $\tilde{K}$ are calculated by the converged values of the solution to the PGF of $\mathbf{P}_2$. }
\label{ADGE}
\end{figure*}
\begin{figure*}[!htbp]
\centering
\includegraphics[width=0.77\linewidth]{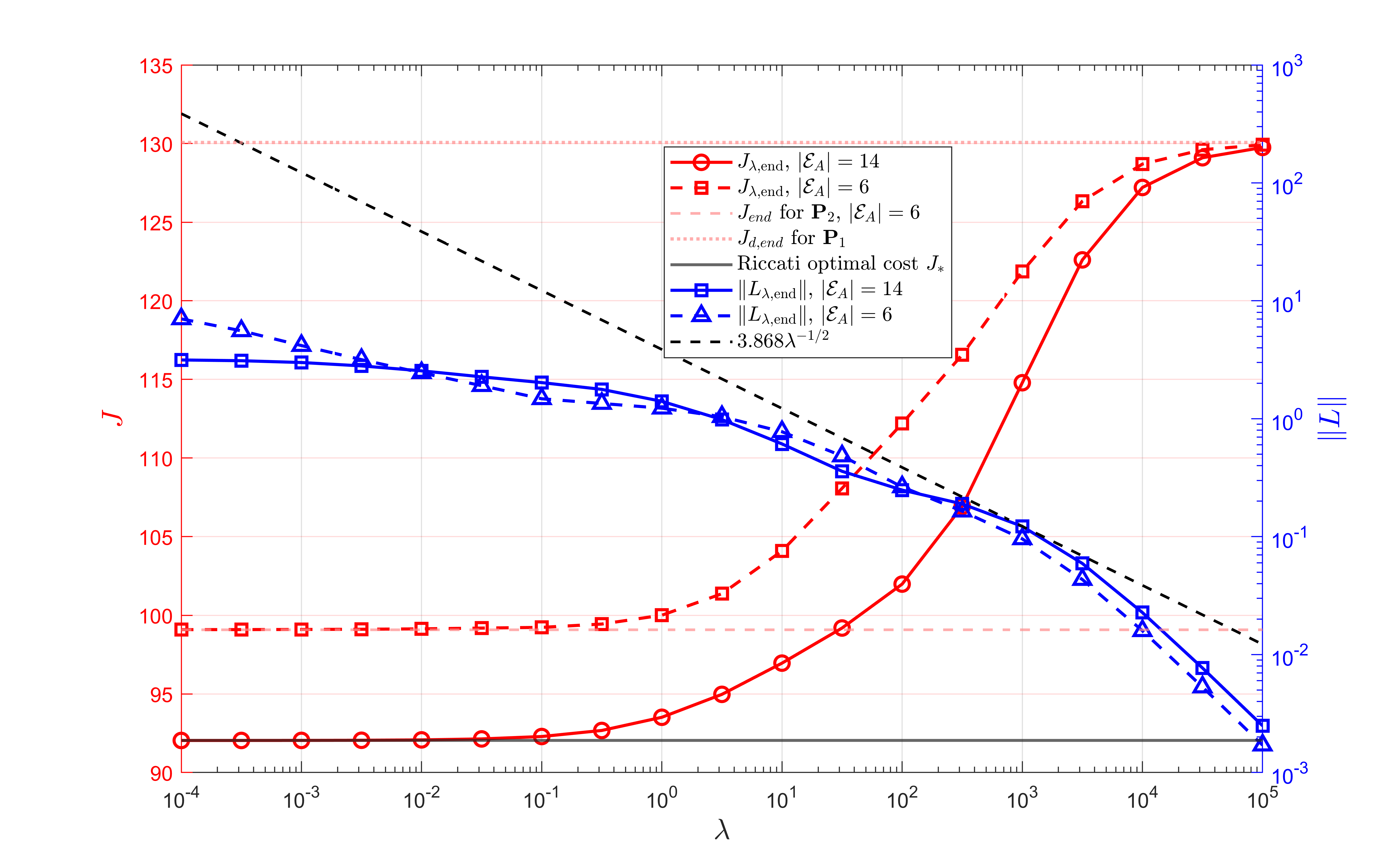}
\caption{Illustration of the effects of penalty parameter $\lambda$ on cost $J_{\lambda}$ and norm $\|L_\lambda\|$. }
\label{REGU}
\end{figure*}

\noindent\textbf{Effects of ADGs on suboptimality gap:} We consider PGFs of $\mathbf{P}_2$ under different ADGs with a common initial policy $K_0$. Although  convergence of the solution to the PGF of $\mathbf{P}_2$ cannot be guaranteed in general, we focus on $\mathbf{P}_2$, since the cost along the PGF always converges, and it is the performance improvement that concerns us here. Thus convergence of the PGF is declared when either the norm of the projected gradient or the variation of the cost drops below $1e-6$. The converged cost is denoted as $J_{end}$. As shown in Figure \ref{ADGE}(a), we choose randomly a sequence of ADGs which finally becomes fully dense. Figure \ref{ADGE}(b) presents the decrease of the suboptimality gap $\frac{J-J_*}{J_*}$ and error $\frac{\|\tilde{K}-\tilde{K}_*\|}{\|\tilde{K}_*\|}$ as the number of edges in $\mathcal{G}_A$ increases from 0 to 15. In fact, a denser ADG provides larger space for the PGF to evolve, and leads to a better performance. But we notice that the cost does not exhibit an obvious decrease as $|\mathcal{E}_A|$ increases from 12 to 15. In particular, there exists no improvement of performance as $|\mathcal{E}_A|$ varies from 14 to 15. The reason is that Assumption \ref{ass5} already holds for $|\mathcal{E}_A|=14$,  and the global optimality is guaranteed. The case of $|\mathcal{E}_A|=15$ corresponds to an auto-regressive policy in MARL, and is  redundant here.

\noindent\textbf{Effects of partial regularization on connecting $\mathbf{P}_1$ and $\mathbf{P}_2$:} Recall that partial regularization  can promote the convergence behavior of the PGF of $\mathbf{P}_2$, as illustrated in the motivating Example 2 and Figure \ref{PGF4}. Building on this observation, we consider the asymptotic behavior of $\mathbf{P}_3$ as $\lambda \to 0$ or $\lambda\to\infty$. Specifically, two cases corresponding to ADGs with $|\mathcal{E}_A|=6$ and $|\mathcal{E}_A|=14$ are discussed. The same initial condition is chosen as before, and the converged values of $J$ and $\|L\|$ are plotted with respect to $\lambda\in[10^{-4},10^{5}]$ in Figure \ref{REGU}.  We observe that in both cases, the converged value $J_{\lambda,end}\to 130.068$, a suboptimal value of $J_d$ as $\lambda \to \infty$. When $|\mathcal{E}_A|=14$,  $J_{\lambda, end}\to J_*$ as $\lambda \to 0$. In contrast, when $|\mathcal{E}_A|=6$,  $J_{\lambda, end}\to 99.083$ as $\lambda\to 0$, a suboptimal value of $\mathbf{P}_2$, which is strictly larger than $J_*$.  Moreover, the norm of  $L_\lambda$ always stays below the blue dashed line $3.868\lambda^{-0.5}$, implying its convergence to 0 as $\lambda \to \infty$. These phenomena are consistent with Theorem \ref{conne}, and reflect the interesting connection between $\mathbf{P}_2$ and $\mathbf{P}_1$. 

\section{Conclusion}

This work introduced ADP to enlarge the policy space and enhance the optimality of distributed policy gradient for multi-agent LQR problems. Two complementary sparsity conditions were proposed to characterize the representability for generic policies. A PGF was then applied to solve the resulting PO problem. To address the possible non-convergence of the PGF caused by the loss of coercivity, a partially regularized problem was introduced and shown to serve as a bridge between the original multi-agent LQR problem and its counterpart with action dependencies. The convergence of the PGF for the regularized problem was established, and a sufficient condition was derived for local convergence to the globally optimal policy. When the complementary sparsity conditions are not satisfied, the ADP-based algorithm has been shown to guarantee performance improvement. However, the extent of this improvement is yet to be precisely quantified. This will be an important direction for our future research.


\appendix

\begingroup
\fontsize{15.5pt}{20pt}\selectfont
\noindent\textbf{\ \ \ \ Appendix}
 \par\endgroup
 
\section{Proofs of Proposition \ref{pro:spars} and Corollary \ref{cor1}}
\label{pro1}
\begin{proof}[Proof of Proposition \ref{pro:spars}] 
1) 
According to the definition of $\mathcal{M}$, $T$ is surjective. Due to the continuity of $T$ and connectivity of $\mathcal{L}\times\mathcal{K}$, $\mathcal{M}=T(\mathcal{L}\times\mathcal{K})$ is also connected. 

2) Let $\tilde{K}\in\mathcal{M}$, if $(L_0,K_0)\in T^{-1}(\tilde{K})$, then $(L_0,tK_0)\in T^{-1}(t\tilde{K})$. Moreover, we have $K_0=(I-L_0)\tilde{K}$. For any $(L,K)\in T^{-1}(\tilde{K})$, it follows that $K=(I-L)\tilde{K}$, and $K-K_0=-(L-L_0)\tilde{K}$, which is a linear matrix equation concerning  $(L-L_0,K-K_0)$. Thus, $T^{-1}(\tilde{K})$ can be viewed as a translation of a linear subspace of $\mathcal{L}\times\mathcal{K}$. In other words, $T^{-1}(\tilde{K})$ is an affine subspace of $\mathcal{L}\times\mathcal{K}$. 

3) According to Lemma \ref{lem:spl}, for any $L\in\mathcal{L}$ and $K\in\mathcal{K}$, we have
   $(I-L)^{-1}K=\sum\limits_{k=0}^{N-1}L^k K$.
Since $L\in \Sparse(G_A \otimes \mathbf{1_m 1^\top_m})$, an induction on $k$ provides that $L^k\in \Sparse(G_A^{k} \otimes \mathbf{1_m 1^\top_m})$ by Lemma \ref{lem:sp}. Note that $K\in \Sparse(G_S \otimes \mathbf{1_m 1^\top_n})$. Then $L^kK\in \Sparse(G_A^{k} G_S\otimes \mathbf{1_m 1^\top_n})$. Hence, we have $$\mathcal{M}\subseteq \Sparse(\sum\limits_{k=0}^{N-1}G^{k}_A G_S\otimes \mathbf{1_m 1^\top_n}).$$ 
 If we choose $K=(I-G_A\otimes \mathbf{1_m 1^\top_m})^{-1}(G_S\otimes \mathbf{1_m 1^\top_n})$, the second statement can be verified. 
\end{proof}

\begin{proof}[Proof of Corollary \ref{cor1}]
    1) According to Proposition \ref{prop:sp}, $\mathcal{M}=\mathcal{K}$ holds if and only if $\Struct (\sum\limits_{k=0}^{N-1}G^{k}_AG_S) = G_S$. Note that $\Struct (G_AG_S )\leq \Struct (\sum\limits_{k=0}^{N-1}G^{k}_AG_S) = G_S$. On the other hand, if $\Struct ( G_A G_S)\leq G_S$, by Lemma \ref{lem:sp}, we know that $$\Struct ( G_A^2 G_S)\leq \Struct (G_A ) \Struct(G_A G_S ) \leq G_A G_S .$$ Then $\Struct(G_A G_S )\leq G_S$ implies $\Struct (G_A^2 G_S )\leq G_S$. By induction, we have  $\Struct (G_A^k G_S )\leq G_S$, which indicates that $\mathcal{M}\subseteq \Sparse(G_S)=\mathcal{K}$. Combining with the fact that  $\mathcal{K} \subseteq \mathcal{M}$, we obtain $\mathcal{K}= \mathcal{M}$.

    2) This can be obtained by result 3) in Proposition \ref{prop:sp}.

    3) The sufficiency part is straightforward. To show the necessity,  for any $\tilde{K}\in\mathbb{R}^{mN\times nN}$, there exists a pair of $(L,K)\in\mathcal{L}\times\mathcal{K}$ such that $(I-L)^{-1}K=\tilde{K}$. Equivalently, we have $K=\tilde{K}-L\tilde{K}$. In a component form, $K_{ij}=\tilde{K}_{ij}-\sum\limits_{k<i}L_{ik}\tilde{K}_{kj}$. Take $\tilde{K}_{ij}=\mathbf{1}_m\mathbf{1}^\top_n$, and $\tilde{K}_{kj}=0$ for $k<i$, then $K_{ij}=\tilde{K}_{ij}\neq 0$, which indicates $(i,j)\in\mathcal{E}_S$. Since $i$ and $j$ can be chosen arbitrarily,  $\mathcal{G}_S$ is a complete graph. 
\end{proof}

\section{Proofs of Proposition \ref{pron}, Theorem \ref{thm:c1} and Corollary \ref{corr}}
\label{appdis0}
\begin{proof}[Proof of Proposition \ref{pron}]
 Note that $(I-L)^{-1}K=\tilde{K}$ is equivalent to $K+L\tilde{K}=\tilde{K}$ due to the invertibility of $I-L$. Furthermore, it can be expressed as 
    \begin{equation}\label{eq:ve0}
       \vvec(K)+(\tilde{K}^\top\otimes I_{mN})\vvec(L)=\vvec(\tilde{K}). 
    \end{equation}
    Based on the definition of $\hat{K}$ and $\hat{L}$, there exist $x\in\mathbb{R}^{n_k}$ and $y\in\mathbb{R}^{n_l}$ such that $\vvec(K)=\hat{K} x$ and $\vvec(L)=\hat{L}y$. Then (\ref{eq:ve0}) is equivalent to 
    \begin{equation}\label{eq:ve00}
       \bigg[\hat{K},(\tilde{K}^\top\otimes I_{mN})\hat{L}\bigg] [x,y]^\top = \vvec(\tilde{K}).
    \end{equation}
  Therefore, (\ref{eq:ve00}) admits a solution, i.e., $\tilde{K}\in\mathcal{M}$,  if and only if the rank condition (\ref{rankc}) holds here. Furthermore, the dimension of the fiber $T^{-1}(\tilde{K})$ for $\tilde{K}\in\mathcal{M}$ coincides with that of the solution space of (\ref{eq:ve00}), which is given by $n_l+n_k-\rank \bigg[\hat{K},(\tilde{K}^\top\otimes I_{mN})\hat{L}\bigg]$.
    \end{proof}

\begin{proof}[Proof of Theorem \ref{thm:c1}]
1) Utilizing the row blocked structure in (\ref{eq:pa}) for $\tilde{K}\in\mathbb{R}^{mN\times nN}$, we have 
    \begin{equation}\label{nc01}
        \tilde{K}_i=K_i+\sum_{j\in A_i}L_{ij}\tilde{K}_{j},\forall i=1,\dots,N.
    \end{equation}
    In particular, for $l\in S^c_{i}$, $K_{il}=0$, which yields that 
    \begin{equation} \label{ncc}
 \tilde{K}_i\Big|_{S^c_{i}}=\sum_{j\in A_i}L_{ij}\tilde{K}_{j}\Big|_{S^c_{i}},\forall i=1,\dots,N.
    \end{equation}
  Equivalently, we obtain
      \begin{equation} \label{nc}
L_i\Big|_{A_{i}}\mathcal{Y}_i(\tilde{K})=\tilde{K}_i\Big|_{S^c_{i}}, \forall i=1,\dots,N,
    \end{equation}
where $L_i\Big|_{A_{i}}$ is defined analogously to $\tilde{K}_i\Big|_{S^c_{i}}$ in \eqref{defyz}. Hence (\ref{nc}) admits a solution if and only if  $\tilde{K}_i\Big|_{S^c_{i}}$ lies in the row space of $\mathcal{Y}_i(\tilde{K})$ for $i=1,\dots,N$, i.e.,
\begin{equation}\label{nc1}
    \rank \big(\mathcal{Y}_i(\tilde{K})\big)= \rank \big(\mathcal{Z}_i(\tilde{K}),\forall i=1,\dots,N.
\end{equation}
Thus we have shown that $\mathcal{M}$ is included in the right hand side (RHS) of the equation in 1). On the other hand, if (\ref{nc1}) is valid, i.e., $\tilde{K}$ lies in the RHS of 1),   then the above discussion indicates that (\ref{ncc}) admits a solution. Namely,  $L_{ij}$ for all $j\in\mathcal{N}_{in,A}(i)$ can be determined, and then for $l\in S_{i}$, $K_{il}$ can be derived from (\ref{nc01}) as 
  \begin{equation}
        K_{il}=\tilde{K}_{il}-\sum_{j\in A_i}L_{ij}\tilde{K}_{jl},\forall i=1,\dots,N.
    \end{equation}
Hence the RHS of 1) is also included in $\mathcal{M}$. This completes the proof of  1).
    
2) We show that Assumption \ref{ass:den} is necessary and sufficient for $\mathbb{R}^{mN\times nN}\setminus S_0$ to be dense and open in $\mathbb{R}^{mN\times nN}$.

    Sufficiency: Under Assumption \ref{ass:den}, $\mathcal{Y}_i(\tilde{K})\in\mathbb{R}^{m d_i\times n\delta_i}$ is a tall matrix, i.e., $m d_i\geq n\delta_i$ . Then we have
    \begin{equation}
    \begin{aligned}
 S^i_{0}&= \{\tilde{K}\in\mathbb{R}^{mN\times nN}:  \rank \big(\mathcal{Y}_i(\tilde{K})\big)< n \delta_i\},\\
 & = \{\tilde{K}\in\mathbb{R}^{mN\times nN}: \text{all}\ n \delta_i\text{-by-}n \delta_i\ \text{submatrices of}\ \mathcal{Y}_i(\tilde{K})\ 
 \text{have zero determinants}\}.
    \end{aligned}
    \end{equation}
 Since each determinant is a nontrivial polynomial in the entries of its submatrix, its zero set has zero measure by \cite[Proposition 0]{zero}, and hence so does their intersection, i.e., $S^i_{0}$. Then as the finite union of $ S^i_{0}$, $S_{0}$ is a closed set with zero measure. Hence $\mathbb{R}^{mN\times nN}\setminus S_{0}$
 is open and dense in $\mathbb{R}^{mN\times nN}$. Moreover, for $\tilde{K}\in\mathbb{R}^{mN\times nN}\setminus S_{0}$, we have $\rank \big(\mathcal{Y}_i(\tilde{K})\big)= n \delta_i$, i.e., $\mathcal{Y}_i(\tilde{K})$ has full column rank for $i=1,\dots,N$. Then (\ref{nc1}) is valid here, which indicates that $\mathbb{R}^{mN\times nN}\setminus S_{0}\subseteq \mathcal{M}$ . 
 
    Necessity: By contradiction, we assume there exists $i_0\in\{1,\dots,N\}$ such that
    \begin{equation}
       m\deg_{in,A}(i_0)+n\deg_{in, S}(i_0)< nN. 
    \end{equation}
 Equivalently, $m d_{i_0}<n\delta_{i_0}$. Then it  naturally holds that $\rank \big(\mathcal{Y}_i(\tilde{K})\big)<n \delta_i$ for arbitrary $\tilde{K}\in\mathbb{R}^{mN\times nN}$.
Namely, $S^{i_0}_0=\mathbb{R}^{mN\times nN}$, and $\mathbb{R}^{mN\times nN}\setminus S_0=\emptyset$, which contradicts the dense property of $\mathbb{R}^{mN\times nN}\setminus S_0$. 

Note that $\mathbb{R}^{mN\times nN}\setminus S_{0}\subseteq \mathcal{M}$ always holds regardless of Assumption \ref{ass:den}. This completes the proof.
\end{proof}

\begin{proof}[Proof of Corollary \ref{corr}] The sufficiency of Assumption \ref{ass:den} directly follows from Theorem \ref{thm:c1}. To show its necessity, we argue by contradiction. Similarly, suppose $i_0\in\{1,\dots,N\}$ satisfies
$m d_{i_0}<n\delta_{i_0}$. From \eqref{nc01}, we have 
\begin{equation}\label{nc111}
    \tilde{K}_{i_0 }=K_{i}+\sum_{j\in A_i}L_{i_0 j}\tilde{K}_{j}.
\end{equation}
    By a dimension count, the number of unknown variables in \eqref{nc111} is  $ m^2\deg_{in,A}(i_0)+mn\deg_{in, S}(i_0)=m^2 d_{i_0}+mn(N-\delta_{i_0})<mnN$ , while \eqref{nc111} imposes  $mnN$ independent constraints. Hence \eqref{nc111} is overdetermined and admits no solution in a generic sense. This contradicts the fact that all policies in $\mathbb{R}^{mN\times nN}$ except
for a set of measure zero are contained in $\mathcal{M}$. Thus the necessity of Assumption \ref{ass:den} is proved.
\end{proof}

\section{Proof of Proposition \ref{pro:dis}}
\label{appdis1}
\begin{proof}
Take $N=2$ and $m=n=1$, two subsystems $x_1(t+1)=2x_1(t)+u_1(t)$, $x_2(t+1)=u_2(t)$, and graphs $\mathcal{G}_A$ and $\mathcal{G}_S$ with
 $\mathcal{E}_A=\{(1,2)\}$ and $\mathcal{E}_S=\{(1,1),(2,2),(2,1)\}$. 
 In a compact form, we have 
$A=\diag(2,0),\ B=I_2$. One can verify that both subsystems are stabilizable.
Furthermore, let
\begin{equation}
    L=\begin{bmatrix}0&0\\ \ell&0\end{bmatrix}\in\mathcal{L},\qquad
K=\begin{bmatrix}k_{11}&k_{12}\\0&k_{22}\end{bmatrix}\in\mathcal{K}.
\end{equation}
Then the closed-loop matrix is
\begin{equation}
   A_{cl}(\ell,k_{11},k_{12},k_{22})
\triangleq A-B(I-L)^{-1}K
=\begin{bmatrix}2-k_{11}&-k_{12}\\[2pt] -\ell k_{11}&-k_{22}-\ell k_{12}\end{bmatrix}.
\end{equation}

Observe that if $k_{11}=0$, the closed-loop matrix is 
upper triangular, and always has an eigenvalue $2$. Hence we have $\rho(A_{cl})\ge 2$, and $\mathcal{F}\subseteq\{(\ell,k_{11},k_{12},k_{22})\in\mathbb{R}^4: k_{11}\neq 0\}$.

\begin{itemize}
  \item When $k_{11}>0$, take $(\ell,k_{11},k_{12},k_{22})=(-0.5,\,1,\,0.5,\,-0.25)$, then $A_{cl}$ and its spectral radius $\rho(A_{cl})$ are 
 $$A_{cl}=\begin{bmatrix}
   1&-0.5\\
   0.5&0.5
   \end{bmatrix},\qquad
  \rho(A_{cl})=\frac{\sqrt3}{2}<1;$$ 

  \item When $k_{11}<0$, take $(\ell,k_{11},k_{12},k_{22})=(6.4,\,-0.5,\,1,\,-5.4)$, then we have
$$A_{cl}  =\begin{bmatrix}2.5&-1\\3.2&-1\end{bmatrix},\qquad
  \rho(A_{cl})\approx0.837<1.$$
\end{itemize}
That is, both $\{(\ell,k_{11},k_{12},k_{22})\in\mathbb{R}^4: k_{11}> 0\}$ and
$\{(\ell,k_{11},k_{12},k_{22})\in\mathbb{R}^4: k_{11}<0\}$ are disjoint open sets, each containing points of
$\mathcal{F}$. Hence $\mathcal{F}$ contains at least two regions separated by the hyperplane $k_{11}=0$, and is disconnected.
\end{proof}

\section{Proof of Proposition \ref{pro:no}}
\label{appdis1.1}
\begin{proof}
Take $N=3$ and $m=n=1$. To show that  $J$ is not coercive on $\mathcal{F}$, let $\mathcal{L}=\{L\in\mathbb{R}^{3\times 3}:L_{ij}=0\ \forall j\geq i\},\ \mathcal{K}=\{K\in\mathbb{R}^{3\times 3}:K_{ij}=0\ \forall j\leq i\}$, and 
    $$A=-\frac{1}{2}I_3,\ B=R=I_3,\ Q= \begin{bmatrix}
3 & 1 & 1\\
1 & 2 & -1\\
1 & -1& 4
\end{bmatrix}.$$ Then, we have $\mathcal{F}=\{(L,K)\in\mathcal{L}\times\mathcal{K}:\rho(\frac{1}{2}I+(I-L)^{-1}K)<1\}$.
Denote $L_k=\begin{bmatrix}
0 & 0 & 0\\
0 & 0 & 0\\
k & 0 & 0
\end{bmatrix}$, and $K=0$. Thus, $(L_k,0)\in\mathcal{F}$. Even if $\|L_k\|\to\infty$ as $k\to\infty$, we always have $J(L_k,0)=J_0(0)<\infty$, which implies that $J$ is not coercive on $\mathcal{F}$. 

To obtain coercivity of $J_\lambda$, take a sequence $\{(L_j,K_j)\}^{\infty}_{j=1}\subseteq \mathcal{F}$, we consider the following cases:
\begin{itemize}
    \item Case 1: $\dist[(L_j,K_j),\partial \mathcal{F}]\to 0$ as $j\to\infty$, where $\partial \mathcal{F}$ denotes the boundary of $\mathcal{F}$. Then $T(L_j,K_j)\to\partial \mathcal{S}$, and the coercivity of $J_0$ implies that $J_\lambda(L_j,K_j)\geq J_0(T(L_j,K_j))=J(L_j,K_j) \to\infty$ as $j\to\infty$.
    \item Case 2: $\|L_j\|+\|K_j\|\to\infty$ as $j\to\infty$. Suppose, on the contrary, there exist $\gamma>0$ and a subsequence $\{(L_j,K_j)\}^{\infty}_{j=1}$ (for simplicity, we use the same notation) such that $J_\lambda(L_j,K_j)\leq \gamma$ for  $j\in\mathbb{N}^*$. From the expression of $J_\lambda$, we have $\|L_j\|^2\leq \frac{2\gamma}{\lambda}$. Then it follows that $\lim_{j\to\infty}\|K_j\|=\infty$. Since the Frobenius norm is sub-multiplicative, we have
    \begin{equation}
    \begin{aligned}
        \|K_j\|\leq \|(I-L_j)^{-1}K_j\|\times\|I-L_j\|
        \leq(1+\sqrt{\frac{2\gamma}{\lambda}})\|T(L_j,K_j)\|,
    \end{aligned}
    \end{equation} 
   which indicates $\|T(L_j,K_j)\|\to\infty$ as $j\to\infty$. Based on the coercivity of $J_0$, we conclude that $J_\lambda(L_j,K_j)\geq J_0(T(L_j,K_j)) \to\infty$ as $j\to\infty$. This is in contradiction with $J_\lambda(L_j,K_j)\leq \gamma$. 
\end{itemize}
\end{proof}

\section{Proofs of Lemma \ref{caldt} and Proposition \ref{pro:grad}}
\label{appdis1.5}
\begin{proof}[Proof of Lemma \ref{caldt}]  Given  $(L,K)\in\mathcal{F}$, the tangent map $DT_{(L,K)}:\mathcal{T}_{(L,K)}\mathcal{F}\to \mathcal{T}_{\tilde{K}}\mathcal{M}$ is calculated as 
 \begin{equation}\label{eq:tan}
 \begin{aligned}
       DT_{(L,K)}(\Delta L,\Delta K)&=\frac{d}{ds}\bigg|_{s=0}T(L+s\Delta L,K+s\Delta K)
       =\frac{d}{ds}\bigg|_{s=0}[I-(L+s\Delta L)]^{-1}(K+s\Delta K)\\
       &=(I-L)^{-1}(\Delta K+\Delta L(I-L)^{-1}K)
       =(I-L)^{-1}(\Delta K+\Delta L\tilde{K}).
 \end{aligned}
 \end{equation}
  Then  the expression of $\ker(DT_{(L,K)})$ in (\ref{dtt}) directly follows from (\ref{eq:tan}). Note that $(I-L)^{-1}$ does not affect $\rank(DT_{(L,K)})$. So we can write the matrix representation of mapping $(\Delta L,\Delta K)\mapsto\Delta K+\Delta L\tilde{K}$ as in the proof of Proposition \ref{pron}, which turns out to be $\bigg[\hat{K},(\tilde{K}^\top\otimes I_{mN})\hat{L}\bigg]$.  
\end{proof}

\begin{proof}[Proof of Proposition \ref{pro:grad}]
1)  Recall that the result in \cite[Proposition 3.8]{bu2019} provides  
    \begin{equation}
        \nabla J_0(\tilde{K})=2[R\tilde{K}-B^\top X(A-B\tilde{K})]Y,
    \end{equation} where $\tilde{K}=(I-L)^{-1}K$, and $Y$ is the unique solution to (\ref{eq:lya}).   Based on the fact $J(L,K)=J_0(T(L,K)$, the chain rule for the calculus of matrix products implies
    \begin{align}
        \frac{\partial J}{\partial L}&=(I-L)^{-\top}\nabla J_0(\tilde{K})K^\top(I-L)^{-\top},\\
        \frac{\partial J}{\partial K}&=(I-L)^{-\top}\nabla J_0(\tilde{K}).
    \end{align}
    Then (\ref{gradf}) indicates the desired expression.
    
    2) The global optimality of $(L_*,K_*)$ follows from the fact that $J(L,K)=J_0(T(L,K))\geq J_0(\tilde{K}_*)=J(L_*,K_*)$ for all $(L,K)\in\mathcal{F}$. Based on  Theorem \ref{thm:c1}, we know $\mathcal{M}\cap \mathcal{B}(\tilde{K},\delta)=\mathcal{B}(\tilde{K},\delta)$ is a manifold for $\tilde{K}\in \mathbb{R}^{mN\times nN}\setminus S_0$. Therefore, $T$ can be viewed as a smooth map between two manifolds near $(L,K)\in T^{-1}(\tilde{K})$. Then the chain rule for the composition map provides (\ref{eq:cha}). Furthermore, (\ref{eq:leb}) is easily derived by leveraging $\nabla J_0(\tilde{K}_*)=0$, $\nabla_\mathcal{F} J(L_*,K_*)=0$, and the Leibniz rule.
\end{proof}

\section{Proof of Lemma \ref{long}}
\label{appdis2}
\begin{proof}
Note that all terms on the RHS of (\ref{eq:gfd}) and (\ref{eq:gfr}) are smooth on $\mathcal{F}$. Based on \cite[Theorem 54]{Sontag}, the existence and uniqueness of the solutions to (\ref{eq:gfd}) and (\ref{eq:gfr}) naturally hold. Furthermore, we assume the maximal existence interval of the solutions to (\ref{eq:gfd}) and (\ref{eq:gfr}) is $[0,\tau_{\max})$. 
We only need to show $\tau_{\max}=\infty$.

For (\ref{eq:gfr}), $J_\lambda$ is coercive, and the sublevel set $J^{-1}_\lambda(-\infty,\gamma]$ is compact according to Proposition \ref{pro:no}. Then $\tau_{\max}=\infty$ follows from \cite[Proposition C.3.6]{Sontag}.

For (\ref{eq:gfd}), we denote $X(t)=[L(t), K(t)]$, and  $\dot{X}=-\nabla_\mathcal{F}J(X)$ for simplicity. Observe that 
\begin{equation}
\begin{aligned}
      \|X(t)-X(0)\|&=\|\int^t_{0}\dot{X}(s)ds\|=\|\int^t_{0}\nabla_\mathcal{F}J(X(s))ds\| \\
      &\leq \int^t_{0}\|\nabla_\mathcal{F}J(X(s))\|ds 
      \leq \Bigg(\int^t_{0}\|\nabla_\mathcal{F}J(X(s))\|^2ds\Bigg)^{\frac{1}{2}}\Bigg(\int^t_{0}1ds\Bigg)^\frac{1}{2}  \\ 
      &=\sqrt{t}\Bigg(J(X_0)-J(X(t))\Bigg)^{\frac{1}{2}},
\end{aligned}
\end{equation}
  where we have used the Cauchy-Schwartz inequality $$\Bigg(\int^b_a f(s)g(s)ds\Bigg)^2\leq \int^b_a f^2(s)ds \int^b_a g^2(s)ds,$$ and $\frac{d}{ds}J(X(s))=-\|\nabla_\mathcal{F}J(X(s))\|^2$ along  (\ref{eq:gfd}). Suppose $\tau_{\max}<\infty$. Then for any finite $t\in [0,\tau_{\max})$, we have
  \begin{equation} \label{eq:bdd}
      \|X(t)\|\leq \|X_0\|+\|X(t)-X(0)\|\leq \|X_0\|+\sqrt{\tau_{\max}} J(X_0).
  \end{equation}
 Thus (\ref{eq:bdd}) indicates that $(L(t), K(t))$ is always bounded for $t\in[0,\tau_{\max})$.  
  On the other hand, 
if along (\ref{eq:gfd}), the closed-loop matrix $A-B\tilde{K}(t)$ approaches the boundary of Schur disk, i.e., $\rho(A-B\tilde{K}(t))\to 1$, the coercivity of $J_0$ implies that $J(L(t),K(t))=J_0(\tilde{K}(t))\to\infty$ as $t\to\tau_{\max}$.  This contradicts to the fact that $J(L_0,K_0)\geq J(L(t),K(t))$ along (\ref{eq:gfd}). Hence the distance between the image $T(L(t), K(t))$  and the boundary of
$\mathcal{F}$ keeps a strictly positive lower bound for $t\in[0,\tau_{\max})$.
Combining with the boundedness of $(L(t), K(t))$, $(L(t), K(t))$ will stay in a compact subset of $\mathcal{F}$ for any $t\in[0,\tau_{\max})$. According to \cite[Proposition C.3.6]{Sontag}, this contradicts the assumption that $[0,\tau_{\max} )$ is the maximal existence interval. Therefore,  $\tau_{\max}=\infty$ holds.
\end{proof}

\section{Proof of Theorem \ref{thm:reg}}
\label{appdis4}
\begin{proof}
   From Proposition \ref{pro:no}, the sublevel set $\Omega_\gamma=J^{-1}_\lambda(-\infty,\gamma]$ is compact for any $\gamma\in\mathbb{R}$. Therefore, $(L(t),K(t))$ is contained in a compact subset $\Omega_{\gamma_0}$ of $\mathcal{F}$ with $\gamma_0=J_\lambda(L_0,K_0)$. According to \cite[Theorem 2.2]{2005gra}, the analyticity of $J_\lambda$ combined with boundedness of $(L(t),K(t))$ implies the existence of a stationary point $(L_\ast,K_\ast)\in\mathcal{F}$ such that the solution $(L(t),K(t))$ along the gradient flow of $J_\lambda$ converges to $(L_\ast,K_\ast)$.

   According to \cite[Łojasiewicz gradient inequality]{Lojasiewicz1965}, there exist a Łojasiewicz exponent $\mu\in [\frac{1}{2},1)$\footnote{Łojasiewicz gradient inequality was first
proved in \cite{Lojasiewicz1965}, and later many researchers try to give a simplified proof. In \cite{paul2019}, this inequality was shown to hold for some $\mu\in[\frac{1}{2},1)$ if $x_*$ is a stationary point of an analytic function.}, constant $c>0$, and a neighborhood $\mathcal{U}\subseteq\mathcal{F}$ of $(L_\ast,K_\ast)$ such that 
   \begin{equation}\label{eq:llk}
       \|\nabla_{\mathcal{F}} J_\lambda(L,K)\|\geq c |J_\lambda(L,K)-J_\lambda(L_\ast,K_\ast)|^\mu
   \end{equation}
   holds for $(L,K)\in \mathcal{U} $. Along PGF of $J_\lambda$, denote
 $X(t)=(L(t), K(t))$, and $X_*=(L_*, K_*)$. Since $\lim_{t\to+\infty} X(t)=X_*$,  $X(t)\in\mathcal{U}$ holds for $t$ sufficiently large. Then we have 
   \begin{equation}\label{eq:ins}
      \frac{d}{dt} \Bigg( J_\lambda(X(t))-J_\lambda(X_\ast)\Bigg)=-\|\nabla_{\mathcal{F}}J_\lambda(X(t))\|^2
      \leq -c^2 \Bigg(J_\lambda(X(t))-J_\lambda(X_\ast)\Bigg)^{2\mu}
   \end{equation}
   for $t$ sufficiently large, where we have used (\ref{eq:llk}) and the fact $J_\lambda(X(t))\geq J_\lambda(X_\ast)$. For simplicity, we denote $\psi(t)\triangleq J_\lambda(X(t))-J_\lambda(X_\ast)$. Then $\dot{\psi}(t)\leq -c^2\psi^{2\mu}$ follows from \eqref{eq:ins}. Equivalently, we have 
   \begin{equation}\label{eqq1}
         \frac{d}{dt}\psi^{1-2\mu}(t)\geq (2\mu-1)c^2,\  \text{if}\ \mu>\frac{1}{2};\
       \frac{d}{dt} \log(\psi(t))\leq -c^2,\  \text{if}\ \mu=\frac{1}{2}.
   \end{equation}
   Integrating both hand sides of \eqref{eqq1} from $0$ to $t$ yields that
   \begin{equation}\label{reff3}
       \psi(t)\leq \Bigg(\psi^{1-2\mu}(0)+(2\mu-1)c^2t\Bigg)^{-\frac{1}{2\mu-1}},\   \text{if}\ \mu>\frac{1}{2};\  
       \psi(t)\leq \psi(0)\exp(-c^2t),\  \text{if}\ \mu=\frac{1}{2}.
   \end{equation} 
   Thus $\psi(t)=\mathcal{O}(t^{-\frac{1}{2\mu-1}})$ holds for $\mu>\frac{1}{2}$, and $\psi(t)=\mathcal{O}(\exp(-c^2t))$ holds for $\mu=\frac{1}{2}$  as $t\to\infty$. 
   
   On the other hand, from \eqref{eq:ins}, it also holds that
   \begin{equation}
         \frac{d}{dt} \psi^{1-\mu} = (1-\mu)\psi^{-\mu} \langle\nabla_{\mathcal{F}} J_\lambda(X),\dot{X}\rangle
         =-(1-\mu)\psi^{-\mu} \|\nabla_{\mathcal{F}} J_\lambda(X)\| \|\dot{X}\| \leq -c (1-\mu) \|\dot{X}\|.
   \end{equation}
  Then we notice that 
  \begin{equation}
       \|X(t)-X_\ast\|=\|\int^\infty_{t}\dot{X}(s)ds\|\leq \int^\infty_{t}\|\dot{X}(s)\|ds
       \leq \frac{1}{c(1-\mu)}\psi^{1-\mu}(t).
      \end{equation}
 Combining with (\ref{reff3}), $\|X(t)-X_\ast\|=\mathcal{O}(t^{-\frac{1-\mu}{2\mu-1}})$ holds for $\mu>\frac{1}{2}$, and $\|X(t)-X_\ast\|=\mathcal{O}(\exp(-c^2 t/2))$ holds for $\mu=\frac{1}{2}$  as $t\to\infty$. This completes the proof.
\end{proof}

\section{Proof of Theorem \ref{conne}}
 \label{bridge}
We point out that different from Problems $\mathbf{P}_1$ and $\mathbf{P}_3$, the minimal costs $J_*$ may not be achievable under Assumption \ref{ass:st} since the sub-level sets of $J$ are not compact. Therefore, $J_*$ is in fact the infimum of $J$ generally, as shown in Example 2.

\begin{proof}
\textbf{Convergence behaviors of $J_\lambda$ and $(L_\lambda,K_\lambda)$ as $\lambda\to 0$:}  Since $J_*$ is the infimum of $J$, there exists a sequence $\{(L_r,K_r)\}^\infty_{r=1}$ in $\mathcal{F}$ such that $\lim\limits_{r\to\infty}J(L_r,K_r)=J_*$. Naturally, for an arbitrary $\epsilon>0$, there exists some $R_1\in\mathbb{N}^*$ such that 
\begin{equation}
   J_*-\epsilon < J(L_r,K_r) <   J_*+\epsilon\ \text{for}\ r\geq R_1.
\end{equation}
Then the following inequality holds for $r\geq R_1$:
   \begin{multline}\label{ineq}
       J(L_r,K_r)-\epsilon<J_*\leq J(L_\lambda,K_\lambda)
       \leq J_{\lambda,\ast}=J(L_\lambda,K_\lambda)+\frac{\lambda}{2}\|L_\lambda\|^{2} \leq J(L_r,K_r)+\frac{\lambda}{2}\|L_r\|^{2}<J_*+\frac{\lambda}{2}\|L_r\|^{2}+\epsilon.
   \end{multline}  
Consequently, we have
    \begin{equation}
\label{pr1}
       -\epsilon < J(L_\lambda,K_\lambda)-J(L_r,K_r)\leq \frac{\lambda}{2}\bigg(\|L_r\|^{2}-\|L_\lambda\|^{2}\bigg)\ \text{for}\ r\geq R_1,   
    \end{equation}
 and
     \begin{equation}\label{pr11}
       0\leq J(L_\lambda,K_\lambda)-J_*\leq \frac{\lambda}{2}\bigg(\|L_r\|^{2}-\|L_\lambda\|^{2}\bigg)+\epsilon\ \text{for}\ r\geq R_1,
\end{equation}
   In particular, (\ref{pr1}) indicates the boundedness of $L_\lambda$, i.e.,  $\|L_\lambda\|^2\leq   \|L_{R_1}\|^2+\frac{2\epsilon}{\lambda}$. Combining this with (\ref{pr11}), we obtain 
      \begin{equation}\label{bdd0}
      J(L_\lambda,K_\lambda)-J_*\to 0\ \text{as}\ \lambda\to 0.
   \end{equation}
   
   If Assumption \ref{ass5}  further holds, we know that $(\bar{L}_*,\bar{K}_*)$ in (\ref{delmi}) is well-defined, and $J_*=J(\bar{L}_*,\bar{K}_*)$. A similar argument as in (\ref{ineq}) yields $\|L_\lambda\|\leq   \|\bar{L}_*\|$. 
Note that the coercivity of $J_0$ indicates the boundedness of $\|T(L_\lambda,K_\lambda)\|$ with respect to $\lambda\to0$. Utilizing the boundedness of $L_\lambda$, we derive that both $L_\lambda$ and $K_\lambda$ are bounded as $\lambda\to0$. Then we argue by contradiction. Suppose there exist a subsequence $\{(L_{j},K_j)\}^{\infty}_{j=1}$, $L_\infty\in\mathcal{L}$, and $K'_\ast\in\mathcal{M}\cap\mathcal{S}$ such that $L_{j}\to L_\infty$, and $T(L_j,K_{j})\to K'_\ast\neq\tilde{K}_\ast$ as $j\to\infty$. Since $\tilde{K}_\ast$ is the unique global minimizer of the centralized LQR, the following holds
\begin{equation}
    J(L_j,K_j)=J_0(T(L_j,K_j))\to J_0(K'_\ast)>J_0(\tilde{K}_\ast)=J_*\  \text{as}\ \lambda\to 0,
\end{equation} 
which contradicts (\ref{bdd0}). Therefore, it must hold that $T(L_{\lambda},K_{\lambda})\to \tilde{K}_{\ast}$ as $\lambda\to0$. Without loss of generality, we assume (up to a subsequence) $K_j\to K_\infty$ as $j\to\infty$, and then $T(L_\infty, K_\infty)=\tilde{K}_{\ast}$. Since $\|L_j\|\leq \|\bar{L}_*\|$, it follows that $\|L_\infty\|\leq\|\bar{L}_*\|$. Since $\bar{L}_*$ is the unique element in the fiber with a minimal $L$-component norm, we have $L_\infty=\bar{L}_*$. Then $K_\infty=\bar{K}_\ast$ follows from  $T(L_\infty, K_\infty)=\tilde{K}_{\ast}$. Hence we conclude that $\lim_{\lambda\to 0}L_\lambda=\bar{L}_*,\ \lim_{\lambda\to 0}K_\lambda=\bar{K}_*$.

\noindent\textbf{Convergence behaviors of $J_\lambda$ and $(L_\lambda,K_\lambda)$ as $\lambda\to \infty$:} 
Now we consider another inequality
\begin{multline}
\label{bdd2}
J(L_\lambda,K_\lambda)=J_0(T(L_\lambda,K_\lambda))\leq J_{\lambda,\ast}=J(L_\lambda,K_\lambda)+\frac{\lambda}{2}\|L_\lambda\|^{2}
\leq J_\lambda(0,\hat{K}_\ast)=J(0,\hat{K}_\ast)=J_{d,\ast},
\end{multline}
for any $\hat{K}_\ast\in\hat{\mathcal{E}}_\ast$. As a consequence, we have 
\begin{equation}\label{pr2}
\|L_\lambda\|^2\leq \frac{2}{\lambda}\bigg(J_{d,\ast}-J(L_\lambda,K_\lambda)\bigg)\to 0\ \text{as}\ \lambda\to\infty,
\end{equation}
i.e., $\lim_{\lambda\to \infty}L_\lambda=0$. Based on (\ref{bdd2}), we have
\begin{equation}\label{bdd3}
\liminf_{\lambda\to\infty} J_0(T(L_\lambda,K_\lambda))\leq \limsup_{\lambda\to\infty} J_0(T(L_\lambda,K_\lambda))\leq J_{d,\ast}.
\end{equation}
Similarly, since $L_\lambda$ is bounded, the coercivity of $J_0$ yields that $K_\lambda$ is bounded as well. Moreover,  for any (sub)sequence $\{K_j\}^\infty_{j=1}$ of $\{K_\lambda\}$, there exist $K_\infty\in\mathcal{K}$ and a subsequence of $\{K_j\}^\infty_{j=1}$ (we adopt the same notation) such that $K_j\to K_\infty$ as $j\to\infty$. Then 
\begin{equation}
\lim_{j\to\infty}J_0(T(L_j,K_j))=J_0(T(0,K_\infty))=J_0(K_\infty)=J_d(K_\infty)\geq J_{d,\ast},
\end{equation}
 which indicates $\liminf_{\lambda\to\infty} J_0(T(L_\lambda,K_\lambda))\geq J_{d,\ast}$. From (\ref{bdd3}), we have \begin{equation}
    \lim_{\lambda\to\infty} J(L_\lambda,K_\lambda)=\lim_{\lambda\to\infty} J_0(T(L_\lambda,K_\lambda))= J_{d,\ast}.
\end{equation} By (\ref{bdd2}), we conclude that $\lim_{\lambda\to\infty}J_{\lambda,\ast}=J_{d,\ast}$.  Similarly, the result $\lim_{\lambda\to \infty}\dist(K_\lambda,\hat{\mathcal{E}}_\ast)=0$ can be proved by the argument of contradiction, and is omitted here. 
\end{proof}

\section{Proof of Proposition \ref{sub}}
\label{psub}
\begin{proof}
Step 1: We show that there exist neighborhood $\mathcal{U}_1$ of $\mathcal{E}_*$, and neighborhood $\mathcal{W}_1$ near $\tilde{K}_*$ such that $T:\mathcal{U}_1\to \mathcal{W}_1$ is surjective.  
According to Theorem \ref{thm:c1}, let $\mathcal{W}_1=\mathcal{B}(\tilde{K}_*,\delta_1)\subseteq\mathcal{M}$ for some $\delta_1>0$. Furthermore, the stability margin in standard LQR results indicates that there exists $\delta_2>0$ such that $\mathcal{B}(\tilde{K}_*,\delta_2)\subseteq \mathcal{S}$. Denote $\delta=\min\{\delta_1,\delta_2\}$, $\mathcal{W}_1=\mathcal{B}(\tilde{K}_*,\delta)\subseteq\mathcal{M}\cap\mathcal{S}$, and $\mathcal{U}_1=T^{-1}(\mathcal{W}_1)$. Observe that $\mathcal{U}_1$ is an open neighborhood of $\mathcal{E}_*$ due to the continuity of $T$, and $T:\mathcal{U}_1\to \mathcal{W}_1$ is surjective.

Step 2: We show that there exists an neighborhood $\mathcal{U}$ of $\mathcal{E}_\ast$ such that $T$ is a constant rank map on $\mathcal{U}$. Specifically,
\begin{equation}\label{eq:ran}
     \rank(DT_{(L,K)})=\rank\bigg[\hat{K},(\tilde{K}_\ast^\top\otimes I_{mN})\hat{L}\bigg]
 \end{equation}
 holds for $(L,K)\in\mathcal{U}$. 
 
 Given $\mathcal{K}$ and $\mathcal{L}$, the rank on the RHS of (\ref{eq:ran}) can be viewed as a function of $\tilde{K}_\ast$. 
 Due to the lower semi-continuity of the rank function, there exists neighborhood $\mathcal{W}_2$ of $\tilde{K}_\ast$ such that 
 \begin{equation}\label{eq:ran1}
    \rank\bigg[\hat{K},\tilde{K}^\top\otimes I_{mN})\hat{L}\bigg]\geq\rank\bigg[\hat{K},(\tilde{K}_\ast^\top\otimes I_{mN})\hat{L}\bigg]
 \end{equation}
 holds for each $\tilde{K}\in\mathcal{W}_2$. On the other hand, we know that $\max_{\tilde{K}\in\mathbb{R}^{mN\times nN}} \rank\bigg[\hat{K},(\tilde{K}^\top\otimes I_{mN})\hat{L}\bigg]$ is always achieved for generic $\tilde{K}\in\mathbb{R}^{nN\times nN}$\footnote{This can be obtained via an argument by contradiction. Suppose there exists a generic $\tilde{K}$ which does not achieve the maximal rank $r$. Then the determinants of all $r$-by-$r$ sub-matrices of $\tilde{K}$ vanish. Hence elements of $\tilde{K}$ satisfy some nontrivial polynomial equation with rational coefficients, which contradicts the genericity of $\tilde{K}$.}. Thus, for generic $\tilde{K}_\ast$, we have
 \begin{equation}\label{eq:ran2}
    \rank\bigg[\hat{K},(\tilde{K}_\ast^\top\otimes I_{mN})\hat{L}\bigg]=\max_{\tilde{K}\in\mathbb{R}^{mN\times nN}} \rank\bigg[\hat{K},(\tilde{K}^\top\otimes I_{mN})\hat{L}\bigg].
 \end{equation}
Combining (\ref{eq:ran1}) with  (\ref{eq:ran2}), $\rank\bigg[\hat{K},\tilde{K}^\top\otimes I_{mN})\hat{L}\bigg]$ is constant for $\tilde{K}\in\mathcal{W}_2$. Due to the continuity of $T$, let $\mathcal{W}=\mathcal{W}_1\cap\mathcal{W}_2$, then $\mathcal{U}=T^{-1}(\mathcal{W})$ is an open neighborhood of $\mathcal{E}_\ast$. Based on Lemma \ref{caldt}, $\rank(DT_{(L,K)})$ is constant on $\mathcal{U}$. 

Therefore, $T:\mathcal{U}\to \mathcal{W}$ is surjective and has constant rank. Based on the global rank theorem\cite[Theorem 4.14]{John2012}, $T:\mathcal{U}\to \mathcal{W}$ is a submersion. 
%
\end{proof}

\section{Proof of Proposition \ref{lem:MBP}}
\label{appdis5}
\begin{proof}
By Proposition \ref{pro:grad},  for $(L_*,K_*)\in\mathcal{E}_\ast$, we have
\begin{equation}\label{eq:he}
    \nabla_{\mathcal{F}}^2 J(L_\ast,K_\ast)= DT_{(L_\ast,K_\ast)}^\ast\circ \nabla^2_{\tilde{K}} J_0(\tilde{K}_\ast)\circ DT_{(L_\ast,K_\ast)}.
\end{equation}
According to \cite[Proposition 3.10]{bu2019lqr}, $\nabla^2_{\tilde{K}} J_0(\tilde{K}_\ast)$ is positive definite. Then (\ref{eq:he}) indicates that $\nabla_{\mathcal{F}}^2 J(L_\ast,K_\ast)$ is positive semi-definite, and 
\begin{equation}\label{mb}
    \ker \nabla_{\mathcal{F}}^2 J(L_\ast,K_\ast)=\ker DT_{(L_\ast,K_\ast)}.
\end{equation}
Based on the definition of tangent space in \cite[Chapter 3]{John2012}, the following holds
\begin{equation} 
\begin{aligned}
\mathcal{T}_{(L_*,K_*)} \mathcal{E}_\ast &= \left\{\dot{ \sigma}(0) \,\middle|\, \sigma: \mathbb{R} \to \mathcal{E}_\ast \text{ is smooth},\ \sigma(0) = (L_\ast,K_\ast) \right\}\\
&=\left\{ (\dot{L}(0),\dot{K}(0))\in\mathcal{L}\times\mathcal{K} \,:\, K(t)=(I-L(t))\tilde{K}_\ast \text{ for }t\in\mathbb{R},\ (L(0),K(0)) = (L_\ast,K_\ast) \right\}\\
&=\left\{ (\dot{L}(0),\dot{K}(0))\in\mathcal{L}\times\mathcal{K} : \dot{K}(0)=-\dot{L}(0)\tilde{K}_\ast \right\}.
\end{aligned}
\end{equation} 
According to Lemma \ref{caldt}, it follows that $\ker DT_{(L_*,K_*)}=\mathcal{T}_{(L_*,K_*)}\mathcal{E}_\ast$. Combining with (\ref{mb}), we obtain
\begin{equation}\label{mb1}
    \ker \nabla_{\mathcal{F}}^2 J(L_\ast,K_\ast)=\mathcal{T}_{(L_\ast,K_\ast)}\mathcal{E}_\ast.
\end{equation}
Hence $J$ has the MB property at each $(L_\ast,K_\ast)\in\mathcal{E}_\ast$.  
Since $\lambda_\ast$ is the smallest positive eigenvalue of $\nabla^2_{\mathcal{F}}J(L_\ast,K_\ast)$, Definition  \ref{MBde} indicates that $J$ satisfies the $\lambda_\ast$-Morse-Bott property ($\lambda_\ast$-MB) at $(L_\ast,K_\ast)$. Based on the equivalence among  PL condition, QG, and MB property of $C^2$ functions \cite[Proposition 2.3, Corollary 2.11]{2025MP}, we  obtain the result in 2). 
\end{proof}

\section{Proof of Theorem \ref{thm:c3}}
\label{appdist2}
\begin{proof}
According to \cite[Theorem 3]{stab2006}, any $(L_\ast,K_\ast)\in\mathcal{E}_\ast$, as a global minimizer of $J$, is Lyapunov stable. Then for any open neighborhood $\mathcal{U}$ of $(L_\ast,K_\ast)$, there exists a sub-neighborhood $\mathcal{U}_0\subseteq\mathcal{U}$ such that the trajectory of the PGF of $J$ will lie in $\mathcal{U}$, i.e., $(L(t),K(t))\in\mathcal{U}$ as long as the initial condition $(L_0,K_0)\in\mathcal{U}_0$ holds. Then the result in 1) follows easily from \cite[Theorem 2.2]{2005gra}.

Since $J$ is $\lambda_\ast$-Morse-Bott at $(L_\ast,K_\ast)\in\mathcal{E}$, $\mathcal{U}$ can be chosen to guarantee that  both the PL inequality (\ref{muPL}) and the QG property (\ref{muQG}) are valid based on Proposition \ref{lem:MBP}. 
Utilizing the PL inequality (\ref{muPL}), it follows that along the PGF of $J$, the following holds
\begin{equation}\label{rat}
\begin{aligned}
   \frac{d}{dt} \bigg(J(L(t),K(t))-J_\ast\bigg)=-\|\nabla_{\mathcal{F}}J(L(t),K(t))\|^2
   &\leq -2\mu\bigg(J(L(t),K(t))-J_\ast\bigg)
\end{aligned}
\end{equation}
for $(L_0,K_0)\in\mathcal{U}_0$. This is due to the fact that the trajectory satisfies $(L(t),K(t))\in\mathcal{U}$ if $(L_0,K_0)\in\mathcal{U}_0$. Integrating both hand sides provides
\begin{equation}
    J(L(t),K(t))-J_\ast\leq \bigg(J(L_0,K_0)-J_\ast\bigg)\exp(-2\mu t).
\end{equation}
Combining the QG property (\ref{muQG}), along the PGF of $J$, we have
\begin{equation}
        \begin{aligned}\label{ine}
     \frac{\mu^\prime}{2}\dist[(L(t),K(t)),\mathcal{E}_\ast]^2\leq J(L(t),K(t))-J_\ast
     \leq \bigg(J(L_0,K_0)-J_\ast\bigg)\exp(-2\mu t)
    \end{aligned}
    \end{equation}
for $(L_0,K_0)\in\mathcal{U}_0$. This completes the proof.
\end{proof}

\section{Proof of Theorem \ref{thm:c6}}
\label{appthm5}
\begin{proof} Note that the result in 1) follows from (\ref{con}) in Theorem \ref{conne}. We will prove the conclusion in 2) as follows.
    First, we show that under Assumption \ref{ass5}, $ \nabla_{\mathcal{F}}^2 J_\lambda(\bar{L}_*,\bar{K}_*)$ is strictly positive definite as long as $\lambda>0$. Note that $\nabla_{\mathcal{F}}^2 J_\lambda(\bar{L}_*,\bar{K}_*)$ in (\ref{eq:leb}) is always positive semi-definite. Suppose $\Delta E=[\Delta L,\Delta K]\in\mathcal{L}\times\mathcal{K}$ satisfies $  \nabla_{\mathcal{F}}^2 J_\lambda(\bar{L}_*,\bar{K}_*)[\Delta E,\Delta E]=0$, i.e.,
\begin{equation}\label{inee}
  \nabla_{\mathcal{F}}^2 J(\bar{L}_\ast,\bar{K}_\ast)[\Delta E,\Delta E]+\lambda  \mathcal{Q}_{\mathcal{L}}[\Delta E,\Delta E]=0.
\end{equation}
Since both terms on the left hand side of (\ref{inee}) are positive semi-definite, it follows that
\begin{equation}\label{ine21}
  \nabla_{\mathcal{F}}^2 J(\bar{L}_*,\bar{K}_*)[\Delta E,\Delta E]=0,\ \lambda  \mathcal{Q}_{\mathcal{L}}[\Delta E,\Delta E]=0.
\end{equation}
The second equation in (\ref{ine21}) provides $\Delta L=\Pi_{\mathcal{L}}(\Delta E)=0$, while the first equation, combining with Lemma \ref{caldt}, provides 
\begin{equation}
\label{eqhe}
\Delta E\in \ker( \nabla_{\mathcal{F}}^2 J(\bar{L}_*,\bar{K}_*))=\ker DT_{(\bar{L}_*,\bar{K}_*)}
=\bigg\{(\Delta L,\Delta K)\in\mathcal{L}\times\mathcal{K}:\Delta K+\Delta L\tilde{K}_\ast=0\bigg\}.
\end{equation}
Substituting $\Delta L=0$ into (\ref{eqhe}) yields that $\Delta K=0$ and $\Delta E=0$. Hence $ \nabla_{\mathcal{F}}^2 J_\lambda(\bar{L}_*,\bar{K}_*)$ is positive definite. 

Second, due to the smoothness of $J_\lambda$ near $(\bar{L}_*,\bar{K}_*)$, there exists a small neighborhood $\mathcal{U}$ of $(\bar{L}_*,\bar{K}_*)$ such that $ \nabla_{\mathcal{F}}^2 J_\lambda(L,K)\succ 0$  for $(L,K)\in\mathcal{U}$. Based on 1), there exists $\delta>0$ such that $ \nabla_{\mathcal{F}}^2 J_\lambda(L_\lambda,K_\lambda)\succ 0$  when $0<\lambda<\delta$. Thus $(L_\lambda,K_\lambda)$ is a nondegenerate minimizer of $J_\lambda$. Hence  $(L_\lambda,K_\lambda)$ is isolated, and its EAS property follows from \cite[Section 4]{stab2006}.
\end{proof}

\section{Parameter choices in Section \ref{num}}
\label{parae}
Matrices $G_S$, $G_A$, and $Q_0$ used in Section \ref{num} are given below.
\begin{equation}\label{defgr}
G_S = \begin{bmatrix}
1 & 1 & 1 & 1 & 1 & 1\\
1 & 1 & 1 & 0 & 1 & 1\\
1 & 0 & 1 & 0 & 1 & 1\\
0 & 1 & 0 & 1 & 0 & 1\\
1 & 0 & 0 & 0 & 1 & 0\\
0 & 0 & 0 & 1 & 0 & 1
\end{bmatrix},\quad
G_A = \begin{bmatrix}
0 & 0 & 0 & 0 & 0 & 0\\
1 & 0 & 0 & 0 & 0 & 0\\
1 & 1 & 0 & 0 & 0 & 0\\
1 & 1 & 1 & 0 & 0 & 0\\
1 & 1 & 1 & 1 & 0 & 0\\
1 & 0 & 1 & 1 & 1 & 0
\end{bmatrix},\quad
Q_0 = \begin{bmatrix}
4 & -1 & -1 & 0 & -1 & -1\\
-1 & 4 & -1 & 0 & -1 & -1\\
-1 & -1 & 5 & -1 & -1 & 0\\
0 & 0 & -1 & 2 & -1 & 0\\
-1 & -1 & -1 & -1 & 4.1 & 0\\
-1 & -1 & 0 & 0 & 0 & 6
\end{bmatrix}.
\end{equation}


\bibliography{reference}
\end{document}